\documentclass{amsart}
\usepackage{amsmath,amsfonts,latexsym,amssymb}

\newtheorem{theorem}{Theorem}
\newtheorem{lemma}{Lemma}
\newtheorem{corollary}{Corollary}

\theoremstyle{remark}
\newtheorem{remark}{Remark}
\newtheorem*{ack}{Acknowledgements}

\begin{document}

\markboth{Ritabrata Munshi}{Subconvexity for twists of $GL(3)$ $L$-functions - B}
\title[Subconvexity for twists of $GL(3)$ $L$-functions - B]{The circle method and bounds for $L$-functions - IV:\\ Subconvexity for twists of $GL(3)$ $L$-functions - B}

\author{Ritabrata Munshi}   
\address{School of Mathematics, Tata Institute of Fundamental Research, 1 Dr. Homi Bhabha Road, Colaba, Mumbai 400005, India.}     
\email{rmunshi@math.tifr.res.in}
\thanks{The author is supported by SwarnaJayanti Fellowship, 2011-12, DST, Government of India.}

\begin{abstract}
Let $\pi$ be a $SL(3,\mathbb Z)$ Hecke-Maass cusp form satisfying the Ramanujan conjecture and the Selberg-Ramanujan conjecture, and let $\chi$ be a primitive Dirichlet character modulo $M$, which we assume to be prime for simplicity. We will prove the following subconvex bound
$$
L\left(\tfrac{1}{2},\pi\otimes\chi\right)\ll_{\pi,\varepsilon} M^{\frac{3}{4}-\frac{1}{1612}+\varepsilon}.
$$
\end{abstract}


\maketitle


\section{Introduction}
\label{intro}

Let $\pi$ be a  Hecke-Maass cusp form for $SL(3,\mathbb Z)$ of type $(\nu_1,\nu_2)$ (see \cite{DB} and \cite{G}). Let $\lambda(m,n)$ be the normalized (i.e. $\lambda(1,1)=1$) Fourier coefficients of $\pi$. The Langlands parameters $(\alpha_1, \alpha_2, \alpha_3)$ for $\pi$, are given by $\alpha_1=-\nu_1-2\nu_2+1$, $\alpha_2=-\nu_1+\nu_2$ and $\alpha_3=2\nu_1+\nu_2-1$. Let $\chi$ be a primitive Dirichlet character modulo $M$. The $L$-function associated with the twisted form $\pi\otimes\chi$ is given by the Dirichlet series 
$$
L(s,\pi\otimes\chi)=\sum_{n=1}^\infty \lambda(1,n)\chi(n)n^{-s}
$$
in the domain $\sigma=\text{Re}(s)>1$. The $L$-function extends to an entire function and satisfies a functional equation with arithmetic conductor $M^3$. Hence the convexity bound is given by 
$$
L\left(\tfrac{1}{2},\pi\otimes\chi\right)\ll_{\pi,\varepsilon} M^{3/4+\varepsilon}.
$$ 
The subconvexity problem for this $L$-function has been solved recently in several special cases in \cite{B}, \cite{Mu1}, \cite{Mu2}, \cite{Mu3} and, more recently in \cite{Mu4}. In this paper we will prove a very general statement in this direction. However we are going to assume that the form $\pi$ satisfies the following conditions:
\\

\begin{quote}
(R) The Ramanujan conjecture $\lambda(m,n)\ll (mn)^{\varepsilon}$.
\end{quote}

\begin{quote}
(RS) The Ramanujan-Selberg conjecture $\text{Re}(\alpha_i)=0$.\\
\end{quote}

\begin{theorem}
\label{mthm}
Let $\pi$ be a  Hecke-Maass cusp form for $SL(3,\mathbb Z)$ satisfying conjectures (R) and (RS). Let $\chi$ be a primitive Dirichlet character modulo $M$. Suppose $M$ is a prime number. Then we have
$$
L\left(\tfrac{1}{2},\pi\otimes\chi\right)\ll_{\pi,\varepsilon} M^{\frac{3}{4}-\frac{1}{1612}+\varepsilon}.
$$\\
\end{theorem}

The primality assumption on $M$ is more a technical convenience than an essential requirement. A more general statement without this assumption can be proved using the technique introduced in this paper. Our primary goal here is to present the ideas as clearly as possible without trying to prove the most general statement or the best possible exponent. The conditions (R) and (RS) are quite serious and their removal is a technical challenge. Indeed unlike the previous papers in this series \cite{Mu}, \cite{Mu4}, we do not need to use Deligne type bounds for exponential sums. Instead of estimating exponential sums, we will be required to solve a counting problem, which we tackle in an elementary manner (without recourse to exponential sums). 
\\

The subconvexity problem for $L$-function twisted by a Dirichlet character, has been studied extensively in the literature. The first instance of such a result is of course the pioneering work of Burgess \cite{Bu}, whose well-known bound
$$
L\left(\tfrac{1}{2},\chi\right)\ll M^{\frac{3}{16}+\varepsilon}
$$
still remains unsurpassed. In the case of degree two $L$-functions the problem was first tackled by Duke, Friedlander and Iwaniec \cite{DFI-1} using the amplification technique. Their result has been extended (e.g. to the case of general $GL(2)$ automorphic forms) and improved by several authors in the last two decades. Our theorem provides a $GL(3)$ analogue of the main result of \cite{DFI-1}.\\ 

The present work substantially differs from the other papers in the series, and one may rightly argue that the way we detect the (diagonal) equation $n=r$ here can hardly be termed a circle method. We use the Petersson trace formula which gives an expansion of the Kronecker delta symbol in terms of the Fourier coefficients of holomorphic forms and the Kloosterman sums. However the basic set up for the proof of Theorem~\ref{mthm} coincides with that in \cite{Mu}, \cite{Mu4} and \cite{Mu0}. In particular, we use an expansion of the Kronecker delta symbol to separate the oscillation of the $GL(3)$ Fourier coefficents from that of the character. The idea of using the Petersson formula as a substitute of the circle method is also exploited in \cite{Mu5} where we deal with Rankin-Selberg $L$-functions. \\

\begin{remark}
The approach in this paper gives an unconditional subconvexity result for twists of the symmetric square lifts of $SL(2,\mathbb{Z})$ holomorphic forms. Since in this case (R) follows from the work of Deligne, and we also know the exact location of the poles of the gamma factors. The location of the poles are satisfactory for our purpose. 
\end{remark}

\begin{remark}
The theorem in fact holds under the weaker assumptions $\lambda(m,n)\ll (mn)^{\theta_o+\varepsilon}$ and $|\text{Re}(\alpha_i)|\leq \eta_o$, with $\theta_o$ and $\eta_o$ sufficiently small. Since we need these parameters to be very small, far from what one can hope to achieve in recent future, we refrain from writing it down explicitly. A case of special interest corresponds to symmetric square lifts of $SL(2,\mathbb{Z})$ Maass forms. In this case, though not sufficient for our purpose, strong bounds are known from the work of Kim and Sarnak ($\theta_o=\eta_o=7/32$).
\end{remark}

\ack{The author wishes to thank professors Valentin Blomer, Philippe Michel, Peter Sarnak and Matthew Young for their interest in this work.}\\

\bigskip
\section{The set up}

\subsection{Petersson formula to detect equation $n=r$}
We will start by explaining the expansion of the Kronecker symbol that we will use. Let $p$ be a prime number and let $k\equiv 3\bmod{4}$ be a positive integer. Let $\psi$ be a character of $\mathbb{F}_p^\times$ satisfying $\psi(-1)=-1=(-1)^k$. So in particular $\psi$ is primitive modulo $p$. The collection of Hecke forms of level $p$, weight $k$ and nebentypus $\psi$ is denoted by $H_k(p,\psi)$, and they form an orthogonal basis of the space of cusp forms $S_k(p,\psi)$. Let
$$
\omega_f^{-1}= \frac{\Gamma(k-1)}{(4\pi)^{k-1}\|f\|^2}
$$
be the spectral weights. The Petersson formula gives
$$
\sum_{f\in H_k(p,\psi)}\omega_f^{-1}\lambda_f(n)\overline{\lambda_f(r)}=\delta(n,r)+2\pi i\sum_{c=1}^\infty \frac{S_\psi(r,n;cp)}{cp}J_{k-1}\left(\frac{4\pi\sqrt{nr}}{cp}\right).
$$
This gives an expansion of the Kronecker delta $\delta(n,r)$ (which is the indicator function of the diagonal $n=r$) in terms of the Kloosterman sums
$$
S_{\psi}(a,b;c)=\sideset{}{^\star}\sum_{\alpha\bmod{c}}\psi(\alpha)e\left(\frac{\alpha a+\bar{\alpha}b}{c}\right)
$$ 
and the (Hecke normalized) Fourier coefficients $\lambda_f(n)$ of holomorphic forms $f$, if $pk$ is taken to be sufficiently large (so that the space $S_k(p,\psi)$ is non-trivial). \\

Let $P$ be a parameter which shall be chosen optimally later as a power of the modulus $M$. Let 
$$
P^\star =\sum_{\substack{P<p<2P\\p\;\text{prime}}}\;\sum_{\psi\bmod{p}}\left(1-\psi(-1)\right) =\sum_{\substack{P<p<2P\\p\;\text{prime}}}\;\phi(p) \asymp \frac{P^2}{\log P}.
$$\\

\begin{lemma}
\label{pet-for}
For $Pk\gg 1$ (sufficiently large), we have
\begin{align}
\label{circ-meth}
&\delta(n,r)=\frac{1}{P^\star}\:\sum_{\substack{P<p<2P\\p\;\text{prime}}}\;\sum_{\psi\bmod{p}}\left(1-\psi(-1)\right)\sum_{f\in H_k(p,\psi)}\omega_f^{-1}\lambda_f(n)\overline{\lambda_f(r)}\\
\nonumber &-\frac{2\pi i}{P^\star}\sum_{\substack{P<p<2P\\p\;\text{prime}}}\;\sum_{c=1}^\infty \frac{1}{cp}\sum_{\psi\bmod{p}}\left(1-\psi(-1)\right)S_\psi(r,n;cp)J_{k-1}\left(\frac{4\pi\sqrt{nr}}{cp}\right).
\end{align}
\end{lemma}

\bigskip

\subsection{Approximate functional equation}
We will apply the formula from Lemma~\ref{pet-for} directly to the sum  
\begin{align}
\label{the-sum}
S^\star(N)&=\mathop{\sum\sum}_{m,n=1}^\infty \lambda(m,n)\chi(n)W\left(\frac{nm^2}{N}\right)V\left(\frac{n}{N}\right)\\
\nonumber &=\sum_{r=1}^\infty\mathop{\sum\sum}_{m,n=1}^\infty \lambda(m,n)\chi(r)\delta(n,r)\:W\left(\frac{nm^2}{N}\right)V\left(\frac{r}{N}\right)
\end{align}
where $W$ is a non-negative smooth function supported in $[1,2]$, satisfying $W^{(j)}\ll_j 1$, and $V$ is a smooth function supported in $[M^{-\kappa},4]$, with $V(x)=1$ for $x\in [2M^{-\kappa},2]$, and satisfying $y^jV^{(j)}(y)\ll_j 1$. Also we take $N$ in the range $M^{3/2-\theta}\leq N\ll M^{3/2+\theta}$, with $\theta>0$. The parameter $\theta$ shall be chosen optimally at the end.\\

In the rest of this section we will explain the relevance of the above sum $S^\star(N)$ in the context of the subconvexity problem under focus. First using the definition of $V$ and the bound
\begin{align}
\label{ram-on-av}
\mathop{\sum\sum}_{m^2n\leq x}|\lambda(m,n)|^2\ll x^{1+\varepsilon}
\end{align}
which follows from the Rankin-Selberg theory, we get
\begin{align}
\label{s*}
S^\star(N)=S(N)+ O(NM^{-\kappa/4+\varepsilon}),
\end{align}
where
\begin{align*}
S(N)=\mathop{\sum\sum}_{m,n=1}^\infty \lambda(m,n)\chi(n)W\left(\frac{nm^2}{N}\right).
\end{align*}
Now we will relate the associated Dirichlet series 
$$
\mathop{\sum\sum}_{m,n=1}^\infty \lambda(m,n)\chi(n)(nm^2)^{-s}
$$
to $L(s,\pi\otimes\chi)$. The series is given by the Euler product
$$
\prod_{p\;\text{prime}}\mathop{\sum\sum}_{u,v=0}^\infty \lambda(p^u,p^v)\chi(p)^vp^{-(2u+v)s}.
$$
For $u,v \geq 1$ we have (the Hecke relations)
$$
\lambda(p^u,p^v)=\lambda(p^u,1)\lambda(1,p^v)-\lambda(p^{u-1},1)\lambda(1,p^{v-1}).
$$
So it follows that
\begin{align*}
&\mathop{\sum\sum}_{u,v=0}^\infty \lambda(p^u,p^v)\chi(p)^vp^{-(2u+v)s}\\
&=\mathop{\sum}_{u=0}^\infty \lambda(p^u,1)p^{-2us}\mathop{\sum}_{v=0}^\infty \lambda(1,p^v)\chi(p)^vp^{-vs}\;\left\{1-\chi(p)p^{-3s}\right\}.
\end{align*}
Consequently we have
\begin{align}
\label{l-series-relate}
L(3s,\chi)\mathop{\sum\sum}_{m,n=1}^\infty \lambda(m,n)\chi(n)(nm^2)^{-s}=L(s,\pi\otimes\chi)L(2s,\tilde{\pi})
\end{align}
for $\sigma>1$. Here $\tilde{\pi}$ denotes the dual form.
\\

Next we consider the integral
$$
I=\frac{1}{2\pi i}\int_{(2)}\Lambda(\tfrac{1}{2}+s,\pi\otimes\chi)\Lambda(1+2s,\tilde{\pi})X^s\frac{\mathrm{d}s}{s}.
$$
The product of the completed $L$-functions appearing above is given by 
$$
\Lambda(s,\pi\otimes\chi)\Lambda(2s,\tilde{\pi})=M^{3s/2}\gamma(s)L(s,\pi\otimes\chi)L(2s,\tilde{\pi})
$$
for some gamma factor $\gamma(s)$ (which is basically a product of six gamma functions depending on the Langlands parameters of $\pi$ and the parity of $\chi$). We only need the fact that there are no poles of $\gamma(s)$ in the region $\sigma\geq 0$. We move the contour, in the definition of $I$, to $\sigma=-1/2$. The residue at $s=0$ is given by
$$
M^{3/4}\gamma(1/2)L(\tfrac{1}{2},\pi\otimes\chi)L(1,\tilde{\pi}).
$$
For the integral at $\sigma=-1/2$, which is at the edge of the critical strip, we use trivial bounds to get
$$
\frac{1}{2\pi i}\int_{(-1/2)}\Lambda(\tfrac{1}{2}+s,\pi\otimes\chi)\Lambda(1+2s,\tilde{\pi})X^s\frac{\mathrm{d}s}{s}=O\left(M^{3/2+\varepsilon}X^{-1/2}\right).
$$\\

On the other hand from \eqref{l-series-relate} it follows that the initial integral $I$ is given by
$$
M^{3/4}\;\mathop{\sum\sum}_{m,n=1}^\infty \frac{\lambda(m,n)\chi(n)}{\sqrt{m^2n}}\frac{1}{2\pi i}\int_{(2)}\gamma(1/2+s)L(\tfrac{3}{2}+3s,\chi)\left(\frac{M^{3/2}X}{nm^2}\right)^s\frac{\mathrm{d}s}{s}.
$$
We set
$$
\mathcal{V}(y)=\frac{1}{2\pi i}\int_{(2)}\gamma(1/2+s)L(\tfrac{3}{2}+3s,\chi)y^{-s}\frac{\mathrm{d}s}{s}.
$$
For $y\geq M^\varepsilon$, we see that $\mathcal{V}(y)\ll M^{-2013}$ by shifting the contour to the right. For $0<y<M^\varepsilon$ we shift the contour to $\sigma=\varepsilon$. Differentiating within the integral sign we get 
$$
y^j\mathcal{V}^{(j)}(y)\ll_j 1.
$$
It follows that 
\begin{align*}
M^{3/4}\gamma(1/2)L(\tfrac{1}{2},\pi\otimes\chi)L(1,\tilde{\pi})=M^{3/4}&\;\mathop{\sum\sum}_{m,n=1}^\infty \frac{\lambda(m,n)\chi(n)}{\sqrt{m^2n}}\mathcal{V}\left(\frac{nm^2}{XM^{3/2}}\right)\\
&+O\left(M^{3/2+\varepsilon}X^{-1/2}\right).
\end{align*}
Since $L(1,\tilde{\pi})\gg 1$, taking a smooth dyadic subdivision, and picking $X=M^{\theta-\varepsilon}$, we conclude the following.\\

\begin{lemma}
\label{approx-form}
We have
$$
L(\tfrac{1}{2},\pi\otimes\chi)\ll M^\varepsilon\sup_N \:\frac{S(N)}{\sqrt{N}}+ M^{3/4-\theta/2+\varepsilon}
$$
where the supremum is taken over $N$ in the range $M^{3/2-\theta}<N<M^{3/2+\theta}$, and the weight function $W$ (appearing in the sum $S(N)$) is non-negative, smooth, supported in $[1,2]$, satisfying $W^{(j)}\ll_j 1$. 
\end{lemma}

Now using \eqref{s*} we get
$$
L(\tfrac{1}{2},\pi\otimes\chi)\ll M^\varepsilon\sup_N \:\frac{S^\star(N)}{\sqrt{N}}+ M^{3/4-\theta/2+\varepsilon}+M^{3/4+(2\theta-\kappa)/4+\varepsilon}.
$$
To match the two error terms we take 
$
\kappa=4\theta.
$\\

\begin{corollary}
\label{first-cor}
We have
$$
L(\tfrac{1}{2},\pi\otimes\chi)\ll M^\varepsilon\sup_N \:\frac{S^\star(N)}{\sqrt{N}}+ M^{3/4-\theta/2+\varepsilon},
$$
where the supremum is taken over $N$ in the range $M^{3/2-\theta}<N<M^{3/2+\theta}$.
\end{corollary}

\bigskip

\subsection{Conclusion}
We have reduced the subconvexity problem to that of obtaining sufficient bounds for the sum $S^\star(N)$ as given in \eqref{the-sum}. Applying \eqref{circ-meth} from Lemma~\ref{pet-for} to \eqref{the-sum} we get two terms, namely
$$
S^\star(N)=\mathcal{F}-2\pi i\:\mathcal{O}
$$
where 
\begin{align}
\label{dual-side}
\mathcal{F}=\frac{1}{P^\star}&\sum_{\substack{P<p<2P\\p\;\text{prime}}}\;\sum_{\psi\bmod{p}}\left(1-\psi(-1)\right)\sum_{f\in H_k(p,\psi)}\omega_f^{-1}\\
\nonumber & \times \mathop{\sum\sum}_{m,n=1}^\infty \lambda(m,n)\lambda_f(n)W\left(\frac{nm^2}{N}\right)\:\sum_{r=1}^\infty \overline{\lambda_f(r)}\chi(r)\:V\left(\frac{r}{N}\right),
\end{align}
and 
\begin{align}
\label{off-diag-1}
\mathcal{O}=&\frac{1}{P^\star}\sum_{\substack{P<p<2P\\p\;\text{prime}}}\;\sum_{\psi\bmod{p}}\left(1-\psi(-1)\right) \mathop{\sum\sum}_{m,n=1}^\infty \lambda(m,n)W\left(\frac{nm^2}{N}\right)\\
\nonumber &\times \sum_{r=1}^\infty \chi(r)V\left(\frac{r}{N}\right)\sum_{c=1}^\infty \frac{S_\psi(r,n;cp)}{cp}J_{k-1}\left(\frac{4\pi\sqrt{nr}}{cp}\right).
\end{align} 
We pick the weight $k$ to be large, say of size $\varepsilon^{-1}$.
The second sum which we call the off-diagonal can be nicely bounded if $P$ is taken sufficiently large. On the other hand to the first sum we will apply the functional equations followed by the Petersson formula. The resulting diagonal term vanishes, and the off-diagonal term (which we will call dual off-diagonal) can be bounded nicely if $P$ is taken in a suitable range. We will show that there is an optimal choice of $P$ for which both the terms can be bounded satisfactorily.\\


\section{The off-diagonal}
\label{fod}

In this section we will analyse the off-diagonal contribution $\mathcal{O}$ as given in \eqref{off-diag-1}. Suppose we take $P\gg N^{1/2+\varepsilon}$. Since we are picking $k$ very large, of the order $\varepsilon^{-1}$, and 
$$
J_{k-1}(x)\ll x^{k-1},
$$
the contribution from the tail $c>N^{1/2-\varepsilon}$ is negligibly small. In particular the contributing $c$ are necessarily coprime with $p$. We make a dyadic subdivision of the the $c$-sum, and extract the oscillation from the Bessel function. This leads us to the study of the sum
\begin{align*}
\mathcal{O}(m)=\sum_{\substack{P<p<2P\\p\;\text{prime}}}\;\sum_{\psi\bmod{p}}&\left(1-\psi(-1)\right) \mathop{\sum}_{n=1}^\infty \lambda(m,n)\sum_{r=1}^\infty \chi(r)\\
&\times \sum_{\substack{c=1\\(c,p)=1}}^\infty S_\psi(r,n;cp)e\left(\frac{2\sqrt{nr}}{cp}\right)W_0\left(\frac{n}{N_0},\frac{r}{R},\frac{c}{C}\right),
\end{align*}
for any fixed $m\leq \sqrt{N}$, where $W_0$ is a smooth function supported in $[1,2]^3$, with $W_0^{(i,j,k)}(x,y,z)\ll_{i,j,k} 1$ and $N_0=N/m^2$, $NM^{-\kappa}\ll R \ll N$, $C\ll \sqrt{N_0R}M^\varepsilon/P$. 
From a bound for $\mathcal{O}(m)$ we can conclude a bound for $\mathcal{O}$ via the inequality
\begin{align}
\label{off-diag-bd-bc}
\mathcal{O}\ll \frac{M^\varepsilon}{P^2}\sum_{1\leq m\ll \sqrt{N}}\:\sup \frac{\mathcal{O}(m)}{\sqrt{CP}(N_0R)^{1/4}}
\end{align}
where the supremum is taken over all $C$ and $R$ in the above ranges. Technically speaking one should also take supremum over a class of weight functions, but that does not affect the bound. Note that here we are using the decomposition 
$$
J_{k-1}(2\pi x)=e(x)W(x)+e(-x)\bar{W}(x)
$$
with $W(x)\ll x^{-1/2}$ for $x>1$.\\

\subsection{Sum over $\psi$ and reciprocity}
Using the coprimality $(c,p)=1$, we get
\begin{align*}
&\sum_{\psi\bmod{p}}\left(1-\psi(-1)\right) S_\psi(r,n;cp)\\
&= S(\bar{p}r,\bar{p}n;c)\sum_{\psi\bmod{p}}\left(1-\psi(-1)\right) S_\psi(\bar{c} r,\bar{c} n;p)\\
&= \phi(p)S(\bar{p}r,\bar{p}n;c)\;\left(e\left(\frac{\bar{c}(r+n)}{p}\right)-e\left(-\frac{\bar{c}(r+n)}{p}\right)\right).
\end{align*}
For notational simplicity we will only focus on the contribution of the first term to $\mathcal{O}(m)$, which is given by
\begin{align*}
\mathcal{O}_1(m)=&\sum_{\substack{P<p<2P\\p\;\text{prime}}}\;\phi(p) \mathop{\sum}_{n=1}^\infty \lambda(m,n)\sum_{r=1}^\infty \chi(r)\\
&\times \sum_{\substack{c=1\\(c,p)=1}}^\infty S(\bar{p}r,\bar{p}n;c)e\left(\frac{\bar{c}(r+n)}{p}\right)e\left(\frac{2\sqrt{nr}}{cp}\right)W_0\left(\frac{n}{N_0},\frac{r}{R},\frac{c}{C}\right).
\end{align*}
Our next step is a conductor lowering mechanism. This is one of the most vital steps. Similar tricks were also used in the series of papers \cite{Mu1}, \cite{Mu2} and \cite{Mu3}. There a part of the Kloosterman sum could be evaluated as the modulus was powerful. Here the extra average over $\psi$ helps us to evaluate precisely the twisted average value of the Kloosterman sum. This also makes way for the application of the reciprocity relation
$$
e\left(\frac{\bar{c}(r+n)}{p}\right)=e\left(-\frac{\bar{p}(r+n)}{c}\right)e\left(\frac{r+n}{cp}\right).
$$ 
We will push the last oscillatory factor into the weight function. Note that this is only mildly oscillating at the transition range for $c$. We set
\begin{align}
\label{w1-expli}
W_1\left(x,y,z\right)=e\left(\frac{Ry+N_0x}{Cpz}\right)e\left(\frac{2\sqrt{RN_0xy}}{Cpz}\right)W_0\left(x,y,z\right),
\end{align}
and note that
\begin{align*}
\frac{\partial^{j_1+j_2}}{\partial x^{j_1}\partial y^{j_2}}W_1\left(x,y,z\right)\ll_{j_1,j_2} \left(\frac{N_0}{CP}+\frac{\sqrt{RN_0}}{CP}+1\right)^{j_1}\left(\frac{R}{CP}+\frac{\sqrt{RN_0}}{CP}+1\right)^{j_2}.
\end{align*}\\

\begin{lemma}
\label{lem-3}
We have
\begin{align*}
\mathcal{O}_1(m)=\sum_{\substack{P<p<2P\\p\;\text{prime}}}\;&\sum_{\substack{c=1\\(c,p)=1}}^\infty\phi(p) \mathop{\sum}_{n=1}^\infty \lambda(m,n)\\
&\times\sum_{r=1}^\infty \chi(r) S(\bar{p}r,\bar{p}n;c)e\left(-\frac{\bar{p}(r+n)}{c}\right)W_1\left(\frac{n}{N_0},\frac{r}{R},\frac{c}{C}\right),
\end{align*}
where $W_1$ is as given in \eqref{w1-expli}.
\end{lemma}

\bigskip

\subsection{First application of the Poisson summation}
Next we will apply the Poisson summation formula on the sum over $r$. Observe that before the application of the reciprocity relation the modulus for the sum was $cpM$. Using the reciprocity relation we have brought it down to $cM$. Now we consider the sum over $r$ in the expression in Lemma~\ref{lem-3}. Splitting into congruence classes modulo $cM$ we obtain
$$
\sum_{a\bmod{cM}}\chi(a) S(\bar{p}a,\bar{p}n;c)e\left(-\frac{\bar{p}(a+n)}{c}\right)\sum_{r\in\mathbb{Z}} W_1\left(\frac{n}{N_0},\frac{(a+rcM)}{R},\frac{c}{C}\right).
$$
By Poisson summation we now get
\begin{align*}
\sum_{a\bmod{cM}}&\chi(a) S(\bar{p}a,\bar{p}n;c)e\left(-\frac{\bar{p}(a+n)}{c}\right)\sum_{r\in\mathbb{Z}}\int_{\mathbb{R}}W_1\left(\frac{n}{N_0},\frac{(a+xcM)}{R},\frac{c}{C}\right)e(-rx)\mathrm{d}x.
\end{align*}
The change of variables $(a+xcM)/R\mapsto y$, reduces the above sum to
\begin{align*}
\frac{R}{cM}\sum_{r\in\mathbb{Z}}\;\Bigl\{\sum_{a\bmod{cM}}\chi(a) S(\bar{p}a,\bar{p}n;c)&e\left(-\frac{\bar{p}(a+n)}{c}+\frac{ar}{cM}\right)\Bigr\}\\
&\times \int_{\mathbb{R}}W_1\left(\frac{n}{N_0},y,\frac{c}{C}\right)e\left(-\frac{rRy}{cM}\right)\mathrm{d}y.
\end{align*}\\

Since $C\ll N^{1/2}<M$ we have $(c,M)=1$ (since we are assuming that $M$ is a prime and $\theta$ to be small, say $\theta<1/4$). So the character sum splits into a product of two character sums
$$
\sum_{a\bmod{M}}\chi(a)e\left(\frac{a\bar{c}r}{M}\right) \sum_{a\bmod{c}}S(\bar{p}a,\bar{p}n;c)e\left(-\frac{\bar{p}(a+n)}{c}+\frac{a\bar{M}r}{c}\right).
$$
Writing the first sum in terms of the Gauss sum, and opening the Kloosterman sum we get
$$
\varepsilon_\chi \chi(c\bar{r})\sqrt{M} \sideset{}{^\star}\sum_{b\bmod{c}}e\left(\frac{\bar{p}n(\bar{b}-1)}{c}\right)\sum_{a\bmod{c}}e\left(\frac{\bar{p}a(b-1)}{c}+\frac{a\bar{M}r}{c}\right),
$$
where $\varepsilon_\chi$ is the sign of the Gauss sum for $\chi$. Next we execute the sum over $a$ to arrive at
\begin{align}
\label{oss}
\varepsilon_\chi \chi(c\bar{r})\:c\sqrt{M} \:e\left(\frac{(\overline{1-\bar{M}pr}-1)\bar{p}n}{c}\right)=\varepsilon_\chi \chi(c\bar{r})\:c\sqrt{M} \:e\left(\frac{(\overline{M-pr})rn}{c}\right).
\end{align}
In particular this means that the character sum vanishes unless $(pr-M,c)=1$. This is a restriction on the $r$ sum. Next we consider the integral. By repeated integration by parts we have
$$
\int_{\mathbb{R}}W_1\left(\frac{n}{N_0},y,\frac{c}{C}\right)e\left(-\frac{rRy}{cM}\right)\mathrm{d}y\ll_j \left[\left(\frac{R}{CP}+\frac{\sqrt{RN_0}}{CP}+1\right)\frac{CM}{rR}\right]^j.
$$
Hence the integral is negligibly small if 
$$
|r|\gg M^\varepsilon\left(\frac{M}{P}+\frac{M\sqrt{N_0}}{P\sqrt{R}}+\frac{CM}{R}\right).
$$
The second term `essentially' dominates the last term as $C\ll \sqrt{N_0R}M^\varepsilon/P$. We set
$$
H=\left(1+\frac{\sqrt{N_0}}{\sqrt{R}}\right)\frac{M^{1+\varepsilon}}{P}.
$$
Using the explicit form of $W_1$ as given in \eqref{w1-expli}, and the second derivative test for exponential integrals we get
\begin{align}
\label{second-der}
\int_{\mathbb{R}}W_1\left(\frac{n}{N_0},y,\frac{c}{C}\right)e\left(-\frac{rRy}{cM}\right)\mathrm{d}x\ll_j \frac{\sqrt{CP}}{(N_0R)^{1/4}}.
\end{align}
A more elaborate analysis can be carried out using the stationary phase method. It turns out that the contribution of the stationary point nullifies the oscillation coming from the additive character in \eqref{oss}, via reciprocity relation. This can be used if one wants a better exponent in the main result. \\

\begin{lemma} 
\label{o1m}
We have 
$$\mathcal{O}_1(m)\ll |\mathcal{O}_1^\star(m)|+M^{-2013}$$ where
\begin{align*}
\mathcal{O}_1^\star(m)=\frac{R}{\sqrt{M}} &\sum_{\substack{P<p<2P\\p\;\text{prime}}}\;\sum_{\substack{c=1\\(c,p)=1}}^\infty\phi(p) \mathop{\sum}_{n=1}^\infty \lambda(m,n)\\
&\times \sum_{\substack{|r|<H\\(pr-M,c)=1}}\chi(c\bar{r}) \:e\left(\frac{(\overline{M-pr})rn}{c}\right)W_1^\star\left(\frac{n}{N_0},\frac{rR}{cM},\frac{c}{C}\right),
\end{align*}
with
$$
W^\star_1(x,y,z)=\int_{\mathbb{R}}W_1\left(x,u,z\right)e\left(-uy\right)\mathrm{d}u.
$$
\end{lemma}

\bigskip

\subsection{Cauchy inequality and second application of Poisson summation}
Using Cauchy's inequality we get
\begin{align}
\label{bd-after-cauchy}
\mathcal{O}_1^\star(m)\ll \frac{RM^\varepsilon}{\sqrt{M}}\:\sqrt{\Lambda_m}\:\sqrt{\Psi}
\end{align}
where 
$$
\Lambda_m=\sum_{n\leq 10N/m^2}|\lambda(m,n)|^2
$$ 
and
$$
\Psi=\mathop{\sum}_{n=1}^\infty \Bigl|
\sum_{\substack{P<p<2P\\p\;\text{prime}}}\;\mathop{\sum\sum}_{\substack{1\leq c<\infty \\|r|<H\\(c,p(pr-M))=1}}\phi(p)\chi(c\bar{r}) \:e\left(\frac{(\overline{M-pr})rn}{c}\right)W_1^\star\left(\frac{n}{N_0},\frac{rR}{cM},\frac{c}{C}\right)\Bigr|^2.
$$
We open the absolute square and apply the Poisson summation formula on the sum over $n$ after splitting the sum into congruence classes modulo $cc'$. This gives
$$
\Psi=N_0\mathop{\mathop{\sum\sum}_{\substack{P<p,p'<2P\\p,p'\;\text{prime}}}\;\mathop{\sum\sum\sum\sum}_{\substack{1\leq c, c'<\infty \\|r|, |r'|<H\\(c,p(pr-M))=1\\(c',p'(p'r'-M))=1}}\;\sum_{n\in\mathbb Z}}_{\overline{(M-pr)}c'r-\overline{(M-p'r')}cr'\equiv n\bmod{cc'}}\phi(p)\phi(p')\chi(cr'\overline{c'r})\;U(n,r,r',c,c'),
$$
where
$$
U(n,r,r',c,c')=\int_{\mathbb{R}}W_1^\star\left(x,\frac{rR}{cM},\frac{c}{C}\right)\bar{W}_1^\star\left(x,\frac{r'R}{c'M},\frac{c'}{C}\right)e\left(\frac{nN_0x}{cc'}\right)\mathrm{d}x.
$$\\

By repeated integration by parts we have
$$
U(n,r,r',c,c')\ll_j \left[\left(\frac{N_0}{CP}+\frac{\sqrt{RN_0}}{CP}+1\right)\frac{C^2}{nN_0}\right]^j.
$$
Hence the integral is negligibly small if 
$$
|n|\gg M^\varepsilon\left(\frac{C}{P}+\frac{C\sqrt{R}}{P\sqrt{N_0}}+\frac{C^2}{N_0}\right).
$$
Since $C\ll \sqrt{N_0R}M^\varepsilon/P$, and $N_0, R\ll N$, we see that the right hand side is dominated by 
$$
NM^\varepsilon/P^2.
$$
So if $P\gg N^{1/2+\varepsilon}$, then the contribution of the non-zero frequencies $n\neq 0$ is negligibly small. Hence
$$
\Psi\ll N_0P^2\frac{CP}{\sqrt{N_0R}}\mathop{\mathop{\sum\sum}_{\substack{P<p,p'<2P\\p,p'\;\text{prime}}}\;\mathop{\sum\sum\sum\sum}_{\substack{1\leq c, c'\ll C \\1\leq |r|, |r'|<H\\(c,p(pr-M))=1\\(c',p'(p'r'-M))=1}}}_{\overline{(M-pr)}c'r-\overline{(M-p'r')}cr'\equiv 0\bmod{cc'}}1 \:+\: M^{-2013}.
$$
The factor $CP/\sqrt{N_0R}$ comes from the size of the weight function (see \eqref{second-der}).\\

We have reduced the problem to counting the number of solutions of the above congruence. This we can estimate quite easily. Let $d=(c,c')$, we write $c=de$ and $c'=de'$, with $(e,e')=1$. The congruence condition now reduces to 
$$
\overline{(M-pr)}e'r-\overline{(M-p'r')}er'\equiv 0\bmod{dee'}.
$$ 
The coprimality $(e,e')=1$ now forces $e|r$ and $e'|r'$. Accordingly we write $r=es$ and $r'=e's'$. We are now left with the congruence condition
$$
(M-p'e's')s-(M-pes)s'\equiv 0\bmod{d}.
$$
Given $d$ $e$, $e'$, $s$ and $s'$, the number of $p$, $p'$ in the range $[P,2P]$ satisfying the congruence is dominated by 
$$
P^2(ess',d)/d.
$$
It follows that
$$
\Psi\ll N_0P^4\frac{CP}{\sqrt{N_0R}}\sum_{1\leq d\ll C}\mathop{\sum\sum}_{1\leq e,e'\ll C/d}\mathop{\sum\sum}_{\substack{1\leq s\ll H/e\\1\leq s'\ll H/e'}}\frac{(ess',d)}{d}\ll N_0P^4H^2\frac{CP}{\sqrt{N_0R}}M^\varepsilon.
$$\\

\begin{lemma}
We have
\begin{align*}
\mathcal{O}_1^\star(m)\ll \sqrt{\Lambda_m}\:N_0^{1/4}R^{3/4}P^{5/2}H\sqrt{C}M^{-1/2+\varepsilon}.
\end{align*}
\end{lemma}

\bigskip

\subsection{Conclusion}
Substituting the above bound into Lemma~\ref{o1m}, we see that the same bound holds for $\mathcal{O}_1(m)$. The same bound in fact holds for $\mathcal{O}(m)-\mathcal{O}_1(m)$ as well. Substituting this bound in \eqref{off-diag-bd-bc} we conclude that
\begin{align*}
\mathcal{O}&\ll \frac{M^\varepsilon}{P^2}\sum_{1\leq m\ll \sqrt{N}}\:\sup \sqrt{\Lambda_m}\sqrt{R}P^2M^{-1/2}H\\ 
&\ll \frac{M^{1/2+\varepsilon}}{P^{3/2}}\sum_{1\leq m\ll \sqrt{N}}\:\sqrt{\Lambda_m}\:\sup \left(\sqrt{R}+\sqrt{N_0}\right)\sqrt{P}\ll \frac{NM^{1/2+\varepsilon}}{P}.
\end{align*}
In the last inequality we have used \eqref{ram-on-av}. \\

We summarize the content of this section in the following lemma. Note that to prove this lemma we required neither of the conditions (R) or (RS).
\begin{lemma}
\label{lem-31}
Let $\mathcal{O}$ be as defined in \eqref{off-diag-1}. Suppose $P\gg N^{1/2+\varepsilon}$ with $\varepsilon>0$. Then we have
\begin{align*}
\mathcal{O}\ll \sqrt{N}\frac{M^{5/4+\theta/2+\varepsilon}}{P}.
\end{align*}\\ 
\end{lemma}


\section{Functional equations}
\label{fe}

In the rest of the paper we will analyse $\mathcal{F}$ which is given by \eqref{dual-side}. We will first take a smooth dyadic partition of unity to replace the weight function $V(r/N)$ by a bump function, which we will again denote by $V(r/R)$, having a dyadic support. More precisely we take $R$ in the range $NM^{-\kappa}\ll R\ll N$, and assume that $V$ is smooth, supported in $[1,2]$ and satisfying $V^{(j)}(x)\ll_j 1$. Next we apply summation formulas to the sums over $(m,n)$ and $r$. The summation formulas will be derived from the respective functional equations. (For the sum over $r$ one may also use the $GL(2)$ Voronoi summation formula directly.)\\

\subsection{Functional equation for $L(s,f\otimes\chi)$ and related summation formula}
Let
$$
\Lambda(s,\bar{f}\otimes\chi)=\left(\frac{M\sqrt{p}}{2\pi}\right)^s\Gamma\left(s+\frac{k-1}{2}\right) L(s,\bar{f}\otimes\chi)
$$
be the completed $L$-function associated with the twisted form $\bar{f}\otimes\chi$. We have the functional equation
\begin{align}
\label{func-eqn-1}
\Lambda(s,\bar{f}\otimes\chi)=i^k\bar{\psi}(M)\chi(p)\:\frac{g^2_\chi \bar{g}_{\psi}}{M\sqrt{p}}\:\lambda_f(p)\:\Lambda(1-s,f\otimes\bar{\chi}).
\end{align}
Here $g_\chi$ and $g_\psi$ are the associated Gauss sums. We will use this functional equation to derive a summation formula for the sum
$$
S=\sum_{r=1}^\infty \lambda_{\bar{f}}(r)\chi(r)V\left(\frac{r}{R}\right),
$$
where $V$ is smooth, supported in $[1,2]$, satisfying $V^j(x)\ll_j 1$. By Mellin inversion we get
$$
S=\frac{1}{2\pi i}\int_{(2)}\tilde{V}(s)R^sL(s,\bar{f}\otimes\chi)\mathrm{d}s.
$$
Using \eqref{func-eqn-1} we get
\begin{align}
\label{fe1}
S=&i^k\bar{\psi}(M)\chi(p)\:\frac{g^2_\chi \bar{g}_{\psi}}{M\sqrt{p}}\:\lambda_f(p)\:\frac{M\sqrt{p}}{2\pi}\\
\nonumber &\times \frac{1}{2\pi i}\int_{(2)}\tilde{V}(s)\left(\frac{4\pi^2 R}{M^2p}\right)^{s}\frac{\Gamma(1-s+\frac{k-1}{2})}{\Gamma(s+\frac{k-1}{2})} L(1-s,f\otimes\bar{\chi})\mathrm{d}s.
\end{align}
Let $\mathcal{U}=\{(U,\tilde{R})\}$ be a smooth dyadic partition of unity, which consists of pairs $(U,\tilde{R})$ with $U:[1,2]\rightarrow \mathbb{R}_{\geq 0}$ smooth and 
$$
\sum_{(U,\tilde{R})} U\left(\frac{r}{\tilde{R}}\right)=1,\;\;\;\text{for}\;\;\;r\in(0,\infty).
$$
Also the collection is such that the sum is locally finite in the sense that for any given $\ell\in\mathbb{Z}$ there are only finitely many pairs with $\tilde{R}\in [2^\ell,2^{1+\ell}]$. We move the contour in \eqref{fe1} to $-\varepsilon$, and expand the $L$-function into a series and then use a smooth dyadic partition of unity $\mathcal{U}$, as above, to get 
\begin{align*}
S=i^k\bar{\psi}(M)\chi(p)g^2_\chi &\bar{g}_{\psi}\:\lambda_f(p)\:\frac{1}{2\pi}\sum_{\mathcal{U}}\sum_{r=1}^\infty \frac{\lambda_f(r)\bar{\chi}(r)}{r}U\left(\frac{r}{\tilde R}\right)\\
&\times \frac{1}{2\pi i}\int_{(-\varepsilon)}\tilde{V}(s)\left(\frac{4\pi^2 rR}{M^2p}\right)^{s}\frac{\Gamma(1-s+\frac{k-1}{2})}{\Gamma(s+\frac{k-1}{2})} \mathrm{d}s.
\end{align*}
The poles of the integrand are located at
$$
s=\frac{k+1}{2}+\ell,\;\;\;\text{where}\;\;\ell=0,1,2,\dots.
$$
For $\tilde{R}\gg M^{2+\varepsilon} P/R$ we shift the contour to the left, and for $\tilde{R}\ll M^{2-\varepsilon} P/R$ we shift the contour to $k/2$. Since $k$ is large of the size $\varepsilon^{-1}$, we see that the contribution from the above ranges is negligibly small. Let $\mathcal{U}^\star$ be the subset of $\mathcal{U}$ consisting of those pairs $(U, \tilde{R})$ which have $\tilde{R}$  in the range $M^{2-\varepsilon}P/R\ll \tilde{R}\ll M^{2+\varepsilon} P/R$.\\

\begin{lemma}
\label{lem-fe1}
We have
\begin{align*}
\sum_{r=1}^\infty &\lambda_{\bar{f}}(r)\chi(r)V\left(\frac{r}{R}\right)=i^k\bar{\psi}(M)\chi(p)g^2_\chi \bar{g}_{\psi}\:\lambda_f(p)\:\frac{1}{2\pi}\sum_{\mathcal{U}^\star}\sum_{r=1}^\infty \frac{\lambda_f(r)\bar{\chi}(r)}{r}U\left(\frac{r}{\tilde R}\right)\\
&\times \frac{1}{2\pi i}\int_{(0)}\tilde{V}(s)\left(\frac{4\pi^2 rR}{M^2p}\right)^{s}\frac{\Gamma(1-s+\frac{k-1}{2})}{\Gamma(s+\frac{k-1}{2})} \mathrm{d}s+O(M^{-2013}),
\end{align*}
where $\mathcal{U}^\star$ is as above.
\end{lemma}

\bigskip

\subsection{Functional equation for $L(s,\pi\otimes f)$ and related summation formula}
Now we consider $L(s,\pi\otimes f)$, which is given by the Dirichlet series 
$$
\mathop{\sum\sum}_{m,n=1}^\infty \lambda(m,n)\lambda_{f}(n)(m^2n)^{-s}
$$
in the region of absolute convergences $\text{Re}(s)>1$. The $L$-function extends to an entire function. The completed $L$-function is given by
$$
\Lambda(s,\pi\otimes f)=p^{3s/2}\gamma\left(s\right) L(s,\pi\otimes f),
$$
where $\gamma(s)$ is a product of six gamma factors of the type $\Gamma((s+\kappa_j)/2)$. Also each $\kappa_j$ satisfies $\text{Re}(\kappa_j)>k/2-2$ (see \cite{HMQ}). We have the functional equation
\begin{align}
\label{func-eqn-2}
\Lambda(s,\pi\otimes f)=\iota\left(\frac{\bar{g}_{\bar{\psi}}}{\sqrt{p}}\overline{\lambda_f(p)}\right)^3 \Lambda(1-s,\tilde{\pi}\otimes\overline{f}),
\end{align}
where $\iota$ is a root of unity which depends only on the weight $k$ and the Langlands parameters of $\pi$.\\

Consider the sum (which we will again temporarily denote by $S$)
$$
S=\mathop{\sum\sum}_{m,n=1}^\infty \lambda(m,n)\lambda_{f}(n)W\left(\frac{m^2n}{N}\right).
$$
By Mellin inversion we get
$$
S=\frac{1}{2\pi i}\int_{(2)}\tilde{W}(s)N^sL(s,\pi\otimes f)\mathrm{d}s.
$$
Using functional equation \eqref{func-eqn-2} we see that $S$ is given by
$$
\iota\psi(-1)\left(\frac{g_{\psi}}{\sqrt{p}}\right)^3\overline{\lambda_f(p)}^3p^{3/2}\frac{1}{2\pi i}\int_{(2)}\tilde{W}(s)\left(\frac{N}{p^{3/2}}\right)^{s}\frac{\gamma(1-s)}{\gamma(s)} L(1-s,\tilde{\pi}\otimes\bar{f})\mathrm{d}s.
$$
We move the contour to $-\varepsilon$, expand the $L$-function into a series and then use a partition of unity $\mathcal{U}$, as above, to get 
\begin{align*}
\iota\psi(-1)g_{\psi}^3\overline{\lambda_f(p)}^3 &\sum_{\mathcal{U}}\mathop{\sum\sum}_{m,n=1}^\infty \frac{\lambda(n,m)\overline{\lambda_f(n)}}{m^2n}U\left(\frac{m^2n}{\tilde N}\right)\\
&\times \frac{1}{2\pi i}\int_{(-\varepsilon)}\tilde{W}(s)\left(\frac{m^2nN}{p^3}\right)^{s}\frac{\gamma(1-s)}{\gamma(s)} \mathrm{d}s.
\end{align*}
As before by moving contours we can show that for $\tilde{N}$ outside the range 
$$
[P^3M^{-\varepsilon}/N,P^3M^\varepsilon/N],
$$ 
the total contribution is negligible. Let $\mathcal{U}^\dagger$ be the subset of $\mathcal{U}$ consisting of those pairs $(U,\tilde{N})$ which have $\tilde{N}$ in the above range. \\

\begin{lemma}
\label{lem-fe2}
We have
\begin{align*}
\mathop{\sum\sum}_{m,n=1}^\infty &\lambda(m,n)\lambda_{f}(n)W\left(\frac{m^2n}{N}\right)=\iota\psi(-1)g_{\psi}^3\overline{\lambda(p)}^3 \sum_{\mathcal{U}^\dagger}\mathop{\sum\sum}_{m,n=1}^\infty \frac{\lambda(n,m)\overline{\lambda_f(n)}}{m^2n}U\left(\frac{m^2n}{\tilde N}\right)\\
&\times \frac{1}{2\pi i}\int_{(0)}\tilde{W}(s)\left(\frac{m^2nN}{p^3}\right)^{s}\frac{\gamma(1-s)}{\gamma(s)} \mathrm{d}s + O(M^{-2013}),
\end{align*}
where $\mathcal{U}^\dagger$ is as above.
\end{lemma}

\bigskip

\subsection{Application of Petersson formula}
Now we will apply Lemma~\ref{lem-fe1} and Lemma~\ref{lem-fe2} to \eqref{dual-side}. This reduces the analyses of the sum in \eqref{dual-side} to that of sums of the type
\begin{align*}
\frac{RN}{MP^5}\sum_{\substack{P<p<2P\\p\;\text{prime}}}\;\chi(p)&\sum_{\psi\bmod{p}}\left(1-\psi(-1)\right)\bar{\psi}(-M)g_\psi^2\sum_{f\in H_k(p,\psi)}\omega_f^{-1} \overline{\lambda_f(p)}^2\\
\nonumber &\times\mathop{\sum\sum}_{m,n=1}^\infty \lambda(n,m)\overline{\lambda_f(n)}W\left(\frac{nm^2}{\tilde{N}}\right)\sum_{r=1}^\infty \bar{\chi}(r)\lambda_f(r)W\left(\frac{r}{\tilde{R}}\right)
\end{align*}
where 
\begin{align}
\label{range-all}
\frac{N}{M^{\kappa}}\ll R\ll N,\;\;\; \frac{P^3}{NM^{\varepsilon}}\ll \tilde{N}\ll \frac{P^3M^\varepsilon}{N}\;\;\;\text{and}\;\;\; \frac{M^{2}P}{RM^\varepsilon}\ll \tilde{R}\ll \frac{M^{2+\varepsilon} P}{R}.
\end{align} 
The leading factor accounts for the sizes of the denominators appearing on the right hand side of the summation formulas in Lemmas~\ref{lem-fe1} and \ref{lem-fe2}, and also the size of the Gauss sum associated with $\chi$. \\

We apply the Petersson formula. The diagonal term vanishes as the equality $r=np^2$ never holds in the above range. The off-diagonal is given by
\begin{align}
\label{dual-off-diag-ini}
&\mathcal{O}_{\text{dual}}=\frac{RN}{MP^5}\sum_{\substack{P<p<2P\\p\;\text{prime}}}\;\chi(p)\sum_{\psi\bmod{p}}\left(1-\psi(-1)\right)\bar{\psi}(-M)g_\psi^2\\
\nonumber \times &\mathop{\sum\sum}_{m,n=1}^\infty\sum_{r=1}^\infty \bar{\chi}(r) \lambda(n,m)\sum_{c=1}^\infty\frac{S_{\psi}(np^2,r;cp)}{cp}J_{k-1}\left(\frac{4\pi\sqrt{nr}}{c}
\right) W\left(\frac{nm^2}{\tilde{N}}\right)W\left(\frac{r}{\tilde{R}}\right).
\end{align}
Since the weight $k$ is large, the contribution of the tail $c>\sqrt{\tilde{R}\tilde{N}_0}M^\varepsilon$ is negligible. Here we are setting $\tilde{N}_0=\tilde{N}/m^2$. It follows that the terms where $p^2|c$ make a negligible contribution. \\

Now let us consider the case where $p\|c$. We write $c=pc'$. In this case the Kloosterman sum splits as 
$$
S_{\psi}(np^2,r;cp)=S_{\psi}(0,\overline{c'}r;p^2)S(n,\bar{p}^2r;c').
$$
The first term on the right hand side vanishes unless $p|r$, and accordingly we write $r=pr'$. It follows that
\begin{align*}
&\sum_{\psi\bmod{p}}\left(1-\psi(-1)\right)\bar{\psi}(-M)g_\psi^2S_{\psi}(0,\overline{c'}r;p^2)\\
&=p^2\sum_{\psi\bmod{p}}\left(1-\psi(-1)\right)\psi(\overline{Mc'}r')g_\psi \\
&=p^2\phi(p)\left\{e\left(\frac{Mc'\overline{r'}}{p}\right)-e\left(-\frac{Mc'\overline{r'}}{p}\right)\right\}.
\end{align*}
So the contribution of those $c$ for which $p\|c$ is dominated by 
\begin{align*}
\frac{RNM^\varepsilon}{MP^4}\sum_{\substack{P<p<2P\\p\;\text{prime}}}\:\mathop{\sum\sum}_{nm^2\ll \tilde{N}}\sum_{r\ll \tilde{R}/P}  |\lambda(n,m)|\sum_{c\ll \sqrt{\tilde{R}\tilde{N}_0}M^\varepsilon/P}\frac{1}{\sqrt{c}}.
\end{align*}
Trivially executing the sums over $p$, $r$ and $c$ we arrive at
\begin{align*}
\frac{R\tilde{R}NM^\varepsilon}{MP^4}\:\frac{\sqrt{MP}}{(RN)^{1/4}}\:\mathop{\sum\sum}_{nm^2\ll \tilde{N}}\:  |\lambda(n,m)|.
\end{align*}
Using the Cauchy inequality and applying \eqref{ram-on-av}, we get that the above sum is dominated by
\begin{align*}
O\left(\frac{RNM^\varepsilon}{MP^4}\:\frac{\sqrt{MP}}{(RN)^{1/4}}\:\tilde{R}\tilde{N}\right)=O\left(\frac{\sqrt{MP}}{(RN)^{1/4}}\:M^{1+\varepsilon}\right).
\end{align*}
Observe that we have used the Weil bound for the Kloosterman sum modulo $c$. One may  avoid the application of the Weil bound by employing the Voronoi summation formula on the $n$-sum and then evaluating the remaining sums trivially. \\

We conclude that
\begin{align*}
\mathcal{O}_{\text{dual}}=\mathcal{O}_{\text{red dual}}+O\left(M^{3/2+\kappa/4+\varepsilon}\sqrt{\frac{P}{N}}\right),
\end{align*}
where the reduced dual off-diagonal $\mathcal{O}_{\text{red dual}}$ is given by an expression similar to \eqref{dual-off-diag-ini} but with the extra coprimality restriction $(c,p)=1$.\\

\begin{lemma}
\label{lemma-nine}
We have
\begin{align}
\label{now-num}
\mathcal{F}\ll \sup\:|\mathcal{O}_{\mathrm{red\: dual}}|+M^{3/2+\theta+\varepsilon}\sqrt{\frac{P}{N}},
\end{align}
where the supremum is taken over all R, $\tilde{R}$, $\tilde{N}$ in the range \eqref{range-all}, and
\begin{align}
\label{reduced-dual}
&\mathcal{O}_{\mathrm{red\: dual}}=\frac{RN}{MP^5}\sum_{\substack{P<p<2P\\p\;\text{prime}}}\;\chi(p)\sum_{\psi\bmod{p}}\left(1-\psi(-1)\right)\bar{\psi}(-M)g_\psi^2\\
\nonumber \times &\mathop{\sum\sum}_{m,n=1}^\infty\sum_{r=1}^\infty \bar{\chi}(r) \lambda(n,m)\sum_{\substack{c=1\\(c,p)=1}}^\infty\frac{S_{\psi}(np^2,r;cp)}{cp}J_{k-1}\left(\frac{4\pi\sqrt{nr}}{c}
\right) W\left(\frac{nm^2}{\tilde{N}}\right)W\left(\frac{r}{\tilde{R}}\right).
\end{align}
\end{lemma}

\bigskip


\section{Dual off-diagonal away from transition}
\label{dod}

\subsection{Sum over $\psi$} It remains to study \eqref{reduced-dual}. The Kloosterman sum factorizes as
$$
S_{\psi}(np^2,r;cp)=S_{\psi}(0,\overline{c}r;p)S(\bar{p}np^2,\bar{p}r;c)=S_{\psi}(0,\overline{c}r;p)S(n,r;c).
$$
Observe the curious separation of the variables $n$ and $p$, which is a consequence of the fact that we are studying a $GL(d_1)\times GL(d_2)$ Rankin-Selberg convolution with $d_1-d_2=2$ (in our case $d_1=3$ and $d_2=1$). This inbuilt separation of variables will play an important structural role in our analysis of $\mathcal{O}_{\text{dual}}$. Moreover we have
\begin{align*}
&\sum_{\psi\bmod{p}}\left(1-\psi(-1)\right)\bar{\psi}(-M)g_\psi^2S_{\psi}(0,\overline{c}r;p)\\
&=p\sum_{\psi\bmod{p}}\left(1-\psi(-1)\right)\psi(\overline{Mc}r)g_\psi \\
&=p\phi(p)\left\{e\left(\frac{Mc\overline{r}}{p}\right)-e\left(-\frac{Mc\overline{r}}{p}\right)\right\}.
\end{align*}
(The sum vanishes unless $(r,p)=1$.) We will only deal with the first term. \\

We will take a smooth dyadic subdivision of the $c$-sum in $\mathcal{O}_{\text{red dual}}$. In this section we will show that the contribution of any such subdivision which is away from the transition range, which is marked by $C\sim \sqrt{\tilde{N}_0\tilde{R}}$, is satisfactory. For larger values of $C$ the trivial estimation suffices as the size of the Bessel function is small due to the large weight. We will see that for smaller size of $C$, one can get away with a relatively easy estimate. To this end we fix $C, m\geq 1$ and consider
\begin{align}
\label{dual=od-red}
&\mathcal{O}(C,m)=\frac{RN}{MP^5}\sum_{\substack{P<p<2P\\p\;\text{prime}}}\;\phi(p)\chi(p)
\mathop{\sum}_{n=1}^\infty\sum_{\substack{r=1\\(r,p)=1}}^\infty \bar{\chi}(r) \lambda(n,m)\\
\nonumber \times &\sum_{\substack{c=1\\(c,p)=1}}^\infty\frac{S(n,r;c)}{c}e\left(\frac{Mc\overline{r}}{p}\right)J_{k-1}\left(\frac{4\pi\sqrt{nr}}{c}
\right) W\left(\frac{nm^2}{\tilde{N}}\right)W\left(\frac{r}{\tilde{R}}\right)W\left(\frac{c}{C}\right).
\end{align}
For notational simplicity we are using the same weight function $W(.)$. The final bound however is not influenced by this choice.\\

\begin{lemma}
\label{lem10}
We have
$$
\mathcal{O}_{\mathrm{red\: dual}}\ll M^\varepsilon\sum_{m=1}^{10\sqrt{\tilde{N}}}\:\sup\:|\mathcal{O}(C,m)|+M^{-2013}
$$
where the supremum is taken over all $C\ll M^\varepsilon\sqrt{\tilde{N}_0\tilde{R}}$.
\end{lemma}

\bigskip

Taking absolute values we get
\begin{align*}
|\mathcal{O}(C,m)|\ll \frac{RN}{CMP^5}&\:\sum_{r\in\mathbb{Z}}\sum_{C<c\leq 2C}W\left(\frac{r}{\tilde{R}}\right)\Bigl|\sum_{\substack{P<p<2P\\p\;\text{prime}\\(cr,p)=1}}\;\phi(p)\chi(p)e\left(\frac{Mc\overline{r}}{p}\right)\Bigr|\\
\nonumber &\times\Bigl|\mathop{\sum}_{n=1}^\infty \lambda(n,m)S(n,r;c)J_{k-1}\left(\frac{4\pi\sqrt{nr}}{c}\right) W\left(\frac{nm^2}{\tilde{N}}\right)\Bigr|.
\end{align*}
This is the point where we use the separation of the variables noted above. Now applying the Cauchy inequality (and exploiting positivity) we get
\begin{align}
\label{red-dual-cauchy}
\mathcal{O}(C,m)\ll \frac{RN}{CMP^5}\;\sqrt{\Theta_1}\:\sqrt{\Theta_2}
\end{align}
where
\begin{align}
\label{theta1}
\Theta_1=\mathop{\sum\sum}_{c,r\in\mathbb{Z}}U\left(\frac{c}{C},\frac{r}{\tilde{R}}\right)&\Bigl|\sum_{\substack{P<p<2P\\p\;\text{prime}\\(cr,p)=1}}\;\phi(p)\chi(p)e\left(\frac{Mc\overline{r}}{p}\right)\Bigr|^2
\end{align}
and
\begin{align}
\label{theta2}
\Theta_2=\sum_{r\in\mathbb{Z}}\sum_{C<c\leq 2C}W\left(\frac{r}{\tilde{R}}\right)&\Bigl|\mathop{\sum}_{1\leq n <2\tilde{N}} \alpha(n)S(n,r;c)J_{k-1}\left(\frac{4\pi\sqrt{nr}}{c}\right)\Bigr|^2.
\end{align}
Here 
$$
\alpha(n)=\lambda(n,m) W\left(\frac{nm^2}{\tilde{N}}\right),
$$ 
and $U$ is a suitable compactly supported weight function on $(0,\infty)^2$. 
\\

\subsection{Bound for $\Theta_1$}
Opening the absolute square in the sum \eqref{theta1} we arrive at 
\begin{align*}
&\mathop{\sum\sum}_{\substack{P<p,p'<2P\\p,p'\;\text{prime}}}\;\phi(p)\phi(p')\chi(p\bar{p'})\mathop{\sum\sum}_{\substack{c,r\in\mathbb{Z}\\(cr,pp')=1}}e\left(\frac{Mc\overline{r}}{p}-\frac{Mc\overline{r}}{p'}\right)U\left(\frac{c}{C},\frac{r}{\tilde{R}}\right).
\end{align*}
The diagonal $p=p'$ contribution is dominated by $P^3C\tilde{R}$. Also the coprimality condition $(c,pp')=1$ can be removed at a cost of an error term of size $P^3C\tilde{R}$, which is dominated by the diagonal contribution. We will now  apply the Poisson summation formula on the off-diagonal. Breaking into congruence classes modulo $pp'$ we arrive at
\begin{align*}
&\mathop{\sum\sum}_{\substack{P<p,p'<2P\\p\neq p'\;\text{prime}}}\;\phi(p)\phi(p')\chi(p\bar{p'})\mathop{\sum\sum}_{\substack{\gamma,\rho\bmod{pp'}\\(\rho,pp')=1}}e\left(\frac{M\gamma\overline{\rho}}{p}-\frac{M\gamma\overline{\rho}}{p'}\right)\\
&\times \mathop{\sum\sum}_{c,r\in\mathbb{Z}}U\left(\frac{\gamma+cpp'}{C},\frac{\rho+rpp'}{\tilde{R}}\right).
\end{align*}
Then by the Poisson summation (and standard rescaling) we get
\begin{align*}
C\tilde{R}&\mathop{\sum\sum}_{\substack{P<p,p'<2P\\p\neq p'\;\text{prime}}}\;\frac{\phi(p)\phi(p')}{(pp')^2}\chi(p\bar{p'})\mathop{\sum\sum}_{c,r\in\mathbb{Z}}\mathop{\sum\sum}_{\substack{\gamma,\rho\bmod{pp'}\\(\rho,pp')=1}}e\left(\frac{M\gamma\overline{\rho}}{p}-\frac{M\gamma\overline{\rho}}{p'}+\frac{c\gamma+r\rho}{pp'}\right)\\
&\times \int_{\mathbb{R}^2}U\left(x,y\right)e\left(-\frac{Cc}{pp'}x-\frac{\tilde{R}r}{pp'}y\right)\mathrm{d}x\mathrm{d}y.
\end{align*}
The complete character sum over $\gamma$ now yields the relation
$$
M\bar{\rho}(p'-p)+c\equiv 0\bmod{pp'}.
$$
Hence the above sum reduces to 
\begin{align}
\label{111}
C\tilde{R}\mathop{\sum\sum}_{\substack{P<p,p'<2P\\p\neq p'\;\text{prime}}}\;&\frac{\phi(p)\phi(p')}{pp'}\chi(p\bar{p'})\mathop{\sum\sum}_{\substack{c,r\in\mathbb{Z}\\(c,pp')=1}}e\left(-\frac{\bar{c}rM(p'-p)}{pp'}\right)\\
\nonumber &\times \int_{\mathbb{R}^2}U\left(x,y\right)e\left(-\frac{Cc}{pp'}x-\frac{\tilde{R}r}{pp'}y\right)\mathrm{d}x\mathrm{d}y.
\end{align}
The integral is negligibly small if $|r|\gg P^2M^\varepsilon/\tilde{R}$ or if $|c|\gg P^2M^\varepsilon/C$. \\

Let $V(x)$ be a smooth bump function with support contained in $[-10,10]$, and such that $V^{(j)}\ll_j 1$. Set $R^\star=P^2M^\varepsilon/\tilde{R}$, and consider the sum
\begin{align*}
\mathop{\sum}_{r\in\mathbb{Z}}e\left(-\frac{\bar{c}rM(p'-p)}{pp'}\right)e\left(-\frac{\tilde{R}r}{pp'}y\right)V\left(\frac{r}{R^\star}\right).
\end{align*}
Here $y$ is a fixed positive number. Applying reciprocity we reduce the above sum to
\begin{align*}
\mathop{\sum}_{r\in\mathbb{Z}}e\left(\frac{\bar{p}\bar{p}'rM(p'-p)}{c}\right)e\left(-\frac{rM(p'-p)}{cpp'}-\frac{\tilde{R}r}{pp'}y\right)V\left(\frac{r}{R^\star}\right).
\end{align*}
We break the sum into congruence classes modulo $c$ and then apply the Poisson summation formula. This gives (after standard rescaling)
\begin{align*}
\frac{R^\star}{c}\mathop{\sum}_{r\in\mathbb{Z}}\;&\sum_{\rho\bmod{c}}e\left(\frac{\bar{p}\bar{p}'M(p'-p) \rho+r\rho}{c}\right)\\
&\times\int_{\mathbb{R}}V(z)e\left(-\frac{R^\star M(p'-p)}{cpp'}z-\frac{\tilde{R}R^\star y}{pp'}z\right)e\left(-\frac{R^\star r}{c}z\right)\mathrm{d}z,
\end{align*}
which reduces to
\begin{align*}
R^\star\mathop{\sum}_{\substack{r\in\mathbb{Z}\\c|rpp'+M(p'-p)}}\;\int_{\mathbb{R}}V(z)e\left(-\frac{R^\star M(p'-p)}{cpp'}z-\frac{\tilde{R}R^\star y}{pp'}z\right)e\left(-\frac{R^\star r}{c}z\right)\mathrm{d}z.
\end{align*}
By repeated integration by parts we see that the integral is bounded by
$$
\ll_j\left(\left(1+\frac{MP}{c\tilde{R}}\right)\frac{cM^\varepsilon}{R^\star |r|}\right)^j.
$$
Hence the integral is negligibly small if
$$
|r|\gg M^\varepsilon\left(\frac{\tilde{R}}{C}+\frac{M}{P}\right)\asymp \frac{M^{1+\varepsilon}}{P}.
$$
It follows that \eqref{111} is dominated by
\begin{align*}
C\tilde{R}R^\star\mathop{\sum\sum}_{\substack{P<p,p'<2P\\p\neq p'\;\text{prime}}}\;\mathop{\sum}_{|c|\ll P^2M^\varepsilon/C}\mathop{\sum}_{\substack{|r|\ll M^{1+\varepsilon}/P\\c|rpp'+M(p'-p)}}\:1+M^{-2013}.
\end{align*}
Since $rpp'+M(p'-p)$ never vanishes the sum is seen to be bounded by
\begin{align*}
M^{1+\varepsilon} CP\tilde{R}R^\star\ll M^{1+\varepsilon} CP^3.
\end{align*}
Since $P>\sqrt{N}$ and $\theta$ is sufficiently small, we find that $\tilde{R}>M$. So the above term is dominated by the diagonal contribution.

\begin{lemma}
\label{lem-theta1}
Suppose $P>\sqrt{N}$ then we have
\begin{align}
\label{bd-theta1}
\Theta_1\ll P^3C\tilde{R}M^\varepsilon.
\end{align}
\end{lemma}

\bigskip

\subsection{Bound for $\Theta_2$}
Opening the absolute square in the sum \eqref{theta2} we arrive at 
\begin{align}
\label{theta2-ana}
\Theta_2=\sum_{C<c\leq 2C}&\mathop{\sum\sum}_{1\leq n,n'<2\tilde{N}}\;\alpha(n)\bar{\alpha}(n')\sum_{r\in\mathbb{Z}}S(n,r;c)S(n',r;c)\\
\nonumber &\times J_{k-1}\left(\frac{4\pi\sqrt{nr}}{c}\right)J_{k-1}\left(\frac{4\pi\sqrt{n'r}}{c}\right)W\left(\frac{r}{\tilde{R}}\right).
\end{align}
We only need to consider the case where $C\ll M^\varepsilon\sqrt{\tilde{R}\tilde{N}}/m$, as the Bessel function is negligibly small otherwise due to the large weight. For $C$ in this range we apply the Poisson summation formula on $r$ with modulus $c$. Now the Fourier transform 
$$
\int_{\mathbb{R}}J_{k-1}\left(\frac{4\pi\sqrt{n\tilde{R}}}{c}x\right)J_{k-1}\left(\frac{4\pi\sqrt{n'\tilde{R}}}{c}x\right)W\left(x\right)e\left(-\frac{\tilde{R}r}{c}x\right)\mathrm{d}x
$$
is bounded by
$$
\ll_j \left[\left(1+\frac{\sqrt{\tilde{R}\tilde{N}}}{mC}\right)\frac{C}{\tilde{R}r}\right]^j
$$
by repeated integration by parts $j$ times. Since $C\ll M^\varepsilon\sqrt{\tilde{R}\tilde{N}}/m$, it follows that the integral is negligibly small if 
$$
|r|\gg \sqrt{\frac{\tilde{N}}{\tilde{R}}}\:\frac{M^\varepsilon}{m}.
$$
Since we are going to choose $P\ll M^{1-\varepsilon}$ we have $\tilde{R}\gg \tilde{N}M^\varepsilon$, and hence the non-zero frequencies $r\neq 0$ make a negligible contribution. The main contribution comes from the zero frequency which is given by 
\begin{align*}
\tilde{R}\sum_{C<c\leq 2C}\frac{1}{c}&\mathop{\sum\sum}_{1\leq n,n'<2\tilde{N}}\;\alpha(n)\bar{\alpha}(n')\sum_{a\bmod{c}}S(n,a;c)S(n',a;c)\\
\nonumber &\times \int_{\mathbb{R}}J_{k-1}\left(\frac{4\pi\sqrt{n\tilde{R} x}}{c}\right)J_{k-1}\left(\frac{4\pi\sqrt{n'\tilde{R}x}}{c}\right)W\left(x\right)\mathrm{d}x.
\end{align*}
The integral is bounded by
\begin{align*}
\int_{\mathbb{R}}J_{k-1}\left(\frac{4\pi\sqrt{n\tilde{R} x}}{c}\right)J_{k-1}\left(\frac{4\pi\sqrt{n'\tilde{R}x}}{c}\right)W\left(x\right)\mathrm{d}x\ll \frac{C}{\sqrt{\tilde{R}}(nn')^{1/4}}.
\end{align*}
The character sum is given by
\begin{align*}
\sum_{a\bmod{c}}S(n,a;c)S(n',a;c)=c\:\mathfrak{c}_c(n-n')
\end{align*}
where $\mathfrak{c}_u(v)$ is the Ramanujan sum modulo $u$. We obtain the bound
\begin{align*}
\Theta_2\ll C\sqrt{\tilde{R}}\sum_{C<c\leq 2C}&\mathop{\sum\sum}_{1\leq n,n'<2\tilde{N}}\;\frac{|\alpha(n)||\alpha(n')|}{(nn')^{1/4}}|\mathfrak{c}_c(n-n')|.
\end{align*}
The Ramanujan sum can be bounded by the gcd i.e. $\mathfrak{c}_c(n-n')\ll (c,n-n')$. Consequently
\begin{align*}
\sum_{c\sim C}|\mathfrak{c}_c(n-n')|\ll \begin{cases} CM^\varepsilon&\text{if}\;\;n\neq n'\\ C^2M^\varepsilon &\text{otherwise}.\end{cases}
\end{align*}
So it follows that
\begin{align*}
\Theta_2 &\ll C\sqrt{\tilde{R}}
M^\varepsilon\left\{C^2\mathop{\sum}_{n\sim \tilde{N}/m^2}\;\frac{1}{n^{1/2}}+C\left(\mathop{\sum}_{n\sim \tilde{N}/m^2}\frac{1}{n^{1/4}}\right)^2\right\}\\
&\ll C^2\sqrt{\tilde{R}\tilde{N}_0} 
\left\{C+\tilde{N}_0\right\}M^\varepsilon.
\end{align*}
Since by our choice $P<M$, it follows that $C+\tilde{N}_0\ll \sqrt{\tilde{R}\tilde{N}_0}M^\varepsilon$. We conclude the following lemma. Note that we have used (R) in the above estimate.\\

\begin{lemma}
\label{lem-theta2}
Suppose $P<M$ then we have
\begin{align*}
\Theta_2 &\ll C^2\tilde{R}\tilde{N}_0M^\varepsilon.
\end{align*}
\end{lemma}
 

\bigskip

\subsection{Estimate for $\mathcal{O}(C,m)$ for $C$ away from transition range}
Plugging the above bounds for $\Theta_i$ from Lemma~\ref{lem-theta1} and \ref{lem-theta2} into \eqref{red-dual-cauchy}, we conclude 
\begin{align}
\label{bd-satisfactory}
\mathcal{O}(C,m) &\ll 
\frac{RNM^\varepsilon}{CMP^5}\:\sqrt{P^3C\tilde{R}}\:\sqrt{C^2\tilde{R}\tilde{N}_0}\ll M^\varepsilon\frac{M\sqrt{NC}}{mP}.
\end{align}
Recall that we have already noted that the Bessel function in \eqref{dual=od-red} is negligibly small, because of the large weight $k$, if $C\gg \sqrt{\tilde{R}\tilde{N}_0}M^\varepsilon$. So we need to analyse, for any given $m$, the contribution of $C$ in the range $C\ll \sqrt{\tilde{R}\tilde{N}_0}M^\varepsilon$. In this range
\begin{align*}
\mathcal{O}(C,m) &\ll M^\varepsilon\frac{M\sqrt{N}}{mP}(\tilde{R}\tilde{N}_0)^{1/4}\ll \sqrt{N}\frac{M^{3/4+3\theta/2+\varepsilon}}{m^{3/2}}.
\end{align*}\\

\begin{lemma}
\label{previous}
Suppose $N^{1/2+\varepsilon}<P<M^{1-\varepsilon}$ then we have
$$
\sum_{m>M^{4\theta}}\:\sup\:|\mathcal{O}(C,m)|\ll \sqrt{N}M^{3/4-\theta/2+\varepsilon}
$$
where the supremum is taken over all $C\ll M^\varepsilon\sqrt{\tilde{N}_0\tilde{R}}$.
Also we have 
\begin{align*}
\sum_{m=1}^{10\sqrt{\tilde{N}}}\:\sup_{C<P^2/M^{1/2+\theta}}\:|\mathcal{O}(C,m)|\ll\sqrt{N}M^{3/4-\theta/2+\varepsilon}.
\end{align*}
\end{lemma} 

\bigskip

Substituting the bound from Lemma~\ref{previous} into Lemma~\ref{lem10}, we derive the following corollary.
\begin{corollary}
\label{cor-previous-10}
Suppose $N^{1/2+\varepsilon}<P<M^{1-\varepsilon}$ then we have
$$
\mathcal{O}_{\mathrm{red\: dual}}\ll M^\varepsilon\sum_{m\leq M^{4\theta}}\:\sup\:\mathcal{O}(C,m)+\sqrt{N}M^{3/4-\theta/2+\varepsilon}
$$
where the supremum is taken over all $C$ in the range 
\begin{align}
\label{range-for-C}
\frac{P^2}{M^{1/2}}\frac{1}{M^{\theta}}<C<\frac{P^2M^{1+\varepsilon}}{m\sqrt{NR}}\ll \frac{P^2}{mM^{1/2}}M^{3\theta+\varepsilon}.
\end{align}
\end{corollary}

\bigskip

Later we will be applying Poisson summation on the sum over $c$. To this end we wish to get rid of the coprimality condition $(c,p)=1$ in \eqref{dual=od-red}. Consider the sum in \eqref{dual=od-red} but with the condition $p|c$ in place of $(c,p)=1$, i.e.
\begin{align}
\label{dual=od-red-cop}
&\mathcal{O}^\dagger(C,m)=\frac{RN}{MP^5}\sum_{\substack{P<p<2P\\p\;\text{prime}}}\;\phi(p)\chi(p)
\mathop{\sum}_{n=1}^\infty\sum_{\substack{r=1\\(p,r)=1}}^\infty \bar{\chi}(r) \lambda(n,m)\\
\nonumber &\times\sum_{c=1}^\infty\frac{S(n,r;cp)}{cp}J_{k-1}\left(\frac{4\pi\sqrt{nr}}{cp}
\right) W\left(\frac{nm^2}{\tilde{N}}\right)W\left(\frac{r}{\tilde{R}}\right)W\left(\frac{cp}{C}\right).
\end{align}
Taking absolute value we get
\begin{align*}
&\mathcal{O}^\dagger(C,m)\leq \frac{RN}{CMP^4}\sum_{\substack{P<p<2P\\p\;\text{prime}}}\;\sum_{r\sim\tilde{R}} \sum_{c\sim C/p}\\
&\times \Bigl|
\mathop{\sum}_{n=1}^\infty \lambda(n,m)S(n,r;cp)J_{k-1}\left(\frac{4\pi\sqrt{nr}}{cp}
\right) W\left(\frac{nm^2}{\tilde{N}}\right)\Bigr|.
\end{align*}
Using positivity we glue $c$ and $p$ to arrive at
\begin{align*}
&\mathcal{O}^\dagger(C,m)\ll M^\varepsilon\frac{RN}{CMP^4}\;\sum_{r\sim\tilde{R}} \sum_{C<c\leq 4C}\\
&\times \Bigl|
\mathop{\sum}_{n=1}^\infty \lambda(n,m)S(n,r;c)J_{k-1}\left(\frac{4\pi\sqrt{nr}}{c}\right) W\left(\frac{nm^2}{\tilde{N}}\right)\Bigr|.
\end{align*}
Applying Cauchy inequality we get
\begin{align*}
\mathcal{O}^\dagger(C,m)\ll M^\varepsilon\frac{RN}{CMP^4}\:\sqrt{C\tilde{R}}\;\sqrt{\Theta_2} \ll M^\varepsilon\frac{RN}{CMP^5}\:\sqrt{P^2C\tilde{R}}\;\sqrt{\Theta_2}.
\end{align*}
This can be absorbed in the bound given in Lemma~\ref{previous}. Hence we can drop the coprimality condition $(c,p)=1$ from \eqref{dual=od-red} at no extra cost. We set
$$
\mathcal{O}^\star(C,m)=\mathcal{O}(C,m)+\mathcal{O}^\dagger(C,m).
$$\\

\begin{corollary}
\label{cor-previous-11}
Suppose $N^{1/2+\varepsilon}<P<M^{1-\varepsilon}$ then we have
$$
\mathcal{O}_{\mathrm{red\: dual}}\ll M^\varepsilon\sum_{m\leq M^{4\theta}}\:\sup\:|\mathcal{O}^\star(C,m)|+\sqrt{N}M^{3/4-\theta/2+\varepsilon}
$$
where the supremum is taken over all $C$ in the range \eqref{range-for-C}.
\end{corollary}

\section{Wild dual off-diagonal in transition}
\label{dod2}

In the rest of the paper we will analyse the contribution of those $C$ which lie in the range \eqref{range-for-C} for any given $m\leq M^{4\theta}$. 
We consider the sum in \eqref{dual=od-red} without the coprimality condition $(c,p)=1$, namely
\begin{align}
\label{dual-od-red-tran}
\mathcal{O}^\star(C,m)=&\frac{RN}{CMP^5}\sum_{\substack{P<p<2P\\p\;\text{prime}}}\;\phi(p)\chi(p)\mathop{\sum\sum}_{\substack{c,r=1\\(p,r)=1}}^\infty \bar{\chi}(r)\\
\nonumber &\times\mathop{\sum}_{n=1}^\infty\lambda(n,m)S(n,r;c)e\left(\frac{Mc\overline{r}}{p}\right) W\left(\frac{c}{C},\frac{n}{\tilde{N}_0},\frac{r}{\tilde{R}}\right),
\end{align}
where 
$$
W(x,y,z)=J_{k-1}\left(\frac{4\pi\sqrt{\tilde{N}_0\tilde{R}yz}}{Cx}
\right) x^{-1}W(x)W\left(y\right)W\left(z\right).
$$
The single variable function $W$ on the right hand side is as given in \eqref{dual=od-red}. In particular $W(x,y,z)$ is smooth, supported in the box $[1,2]^3$, and satisfies 
$$
W^{(j_1,j_2,j_3)}(x,y,z)\ll_{j_1,j_2,j_3} M^{\hat{\theta}(j_1+j_2+j_3)}, 
$$
where $\hat{\theta}=4\theta$. 
Moreover it is independent of $p$.\\

\subsection{Voronoi summation formula}
The next step involves an application of the Voronoi summation formula (see \cite{L}, \cite{MS}) on the sum over $n$. Let 
$$
\tilde W(x,s,z)=\int_0^\infty W\left(x,y,z\right)y^{s-1}\mathrm{d}y,
$$ 
and for $\ell=0,1$ define
\begin{align}
\label{gammma-factor}
\gamma_{\ell}(s)=\frac{1}{2\pi^{3(s+\frac{1}{2})}}\prod_{i=1}^3\frac{\Gamma\left(\frac{1+s+\alpha_i+\ell}{2}\right)}{\Gamma\left(\frac{-s-\alpha_i+\ell}{2}\right)}
\end{align}
and set $\gamma_\pm(s)=\gamma_0(s)\mp i\gamma_1(s)$. We define the integral transforms 
\begin{align}
\label{gl}
W^\star_\pm (x,y,z)=\frac{1}{2\pi i}\int_{(\sigma)}y^{-s}\gamma_{\pm}(s)\tilde W(x,-s,z)\mathrm{d}s
\end{align}
where $\sigma>-1+\max\{-\text{Re}(\alpha_1),-\text{Re}(\alpha_2),-\text{Re}(\alpha_3)\}$. Using the bounds for the derivatives of $W$, and using integration by parts we get
$$
\tilde W(x,s,z)\ll_j \frac{M^{\hat{\theta}j}}{|s(s+1)\dots(s+j-1)|}.
$$ 
We can now obtain a bound for the integral transform in \eqref{gl} by shifting the contour to the right and using the Stirling approximation. It follows that $W^\star_\pm (x,y,z)$ is negligibly small if $y\gg M^{3\hat{\theta}+\varepsilon}$. For $0<y\ll M^{3\hat{\theta}+\varepsilon}$ we shift the contour to the left upto  $\sigma=-1+\varepsilon$. Since we are assuming (RS) there are no poles of the gamma factor in this domain. Differentiating under the integral sign we get
\begin{align}
\label{gl-bd}
y^j\frac{\partial^j}{\partial y^j} W^\star_\pm (x,y,z)\ll M^{\hat{\theta}(j-1/2)+\varepsilon}
\end{align}
for $j\geq 1$.\\

\begin{lemma}
For $W$ and $W^\star_\pm$ as above, we have
\begin{align*}
\mathop{\sum}_{n=1}^\infty \lambda(n,m) &\:e\left(\frac{\alpha n}{c}\right)W\left(\frac{c}{C},\frac{n}{\tilde{N}_0},\frac{r}{\tilde{R}}\right)\\
=&c\sum_\pm \sum_{m'|cm}\sum_{n=1}^\infty \frac{\lambda(m',n)}{m'n}S(m\bar{\alpha},\pm n; mc/m')W^\star_\pm\left(\frac{c}{C},\frac{m'^2n\tilde{N}_0}{c^3m},\frac{r}{\tilde{R}}\right).
\end{align*}
\end{lemma}

As we observed above, the tail $m'^2n\gg C^3m^3M^{3\hat{\theta}+\varepsilon}/\tilde{N}$ makes a negligible contribution as the integral transform is negligibly small. For smaller values of $m'^2n$ we take a smooth dyadic subdivision of the $n$-sum, and a dyadic subdivision of the sum over $m'$, to arrive at (consider only the term with $+$ sign)
\begin{align}
\label{after-voronoi}
\frac{\tilde{N}_0}{c^2} \sum_{\substack{m'|cm\\m'\sim\mathfrak{m}'}}\frac{m'}{m}\sum_{n=1}^\infty \lambda(m',n)S(m\bar{\alpha}, n; mc/m')V\left(\frac{c}{C},\frac{n}{L},\frac{r}{\tilde{R}}\right).
\end{align}
Here $V$ is smooth, supported in $[1,2]^3$ and it satisfies
$$
V^{(j_1,j_2,j_3)}(x,y,z)\ll_{j_1,j_2,j_3} M^{\hat{\theta}(j_1+j_2+j_3)}. 
$$
Also (because of \eqref{gl-bd}) we just need to take 
\begin{align}
\label{l-range}
1\leq L\ll \frac{C^3m^3M^{3\hat{\theta}+\varepsilon}}{\mathfrak{m}'^2\tilde{N}}.
\end{align}
The function $V$ involves the (latent) variables $m$ and $m'$, but does not depend on $p$. \\

We have applied the Voronoi summation after opening the Kloosterman sum in the initial expression \eqref{dual-od-red-tran}. So we eventually get the Fourier transform of the Kloosterman sum in \eqref{after-voronoi} which is given by 
$$
\sideset{}{^\star}\sum_{\alpha\bmod{c}}e\left(\frac{\bar{\alpha}r}{c}\right)S(m\bar{\alpha},n; mc/m')=\sideset{}{^\star}\sum_{\beta\bmod{mc/m'}}e\left(\frac{\bar{\beta}n}{mc/m'}\right)\sideset{}{^\star}\sum_{\alpha\bmod{c}}e\left(\frac{\bar{\alpha}(r+\beta m')}{c}\right).
$$
The last sum is a Ramanujan sum. Substituting explicit formula for this sum we obtain
$$
c\sum_{d|c}\frac{\mu(d)}{d}\sideset{}{^\star}\sum_{\substack{\beta\bmod{mc/m'}\\r+\beta m'\equiv 0\bmod{c/d}}}e\left( \frac{\bar{\beta}n}{mc/m'}\right).
$$
Substituting \eqref{after-voronoi} in \eqref{dual-od-red-tran} it follows that to obtain bounds for the expression in \eqref{dual-od-red-tran} we now need to analyse sums of the form 
\begin{align}
\label{dual-od-red-tran-2}
&\mathcal{O}(C,m;L,\mathfrak{m}',\mathfrak{d})=\frac{RN}{CMP^5}\sum_{\substack{P<p<2P\\p\;\text{prime}}}\;\phi(p)\chi(p)\sum_{\substack{r=1\\(p,r)=1}}^\infty \bar{\chi}(r)\mathop{\sum}_{c=1}^\infty e\left(\frac{Mc\overline{r}}{p}\right)\\
\nonumber \times \frac{\tilde{N}_0}{c}& \sum_{\substack{m'|cm\\m'\sim\mathfrak{m}'}}\frac{m'}{m}\sum_{\substack{d|c\\d\sim \mathfrak{d}}}\frac{\mu(d)}{d}\sideset{}{^\star}\sum_{\substack{\beta\bmod{mc/m'}\\r+\beta m'\equiv 0\bmod{c/d}}}\;\sum_{n=1}^\infty \lambda(m',n)e\left( \frac{\bar{\beta}n}{mc/m'}\right)V\left(\frac{c}{C},\frac{n}{L},\frac{r}{\tilde{R}}\right).
\end{align}
Here the weight $V$ is as given in \eqref{after-voronoi}. 

\begin{lemma}
\label{lemma15}
We have
\begin{align*}
\mathcal{O}^\star(C,m)\ll M^\varepsilon \sup |\mathcal{O}(C,m;L,\mathfrak{m}',\mathfrak{d})|
\end{align*}
where the supremum is taken over all $L$ in the range \eqref{l-range}, $\mathfrak{m}'\ll Cm$ and $\mathfrak{d}\ll C$.
\end{lemma}

\bigskip

\subsection{Repeating Voronoi summation}
In the rest of this section we will obtain a bound for \eqref{dual-od-red-tran-2}, which will be satisfactory for our purpose when either $\mathfrak{m}'$ or $\mathfrak{d}$ is suitably large. We call these terms wild. Consider the expression in the second line of \eqref{dual-od-red-tran-2}. Suppose we again apply the Voronoi summation formula on the sum over $n$. (This is the standard reversal process to get rid of the `wild' terms.) Then we arrive at 
\begin{align}
\label{again-voronoi}
\tilde{N}_0&\sum_{\substack{d|c\\d\sim \mathfrak{d}}}\frac{\mu(d)}{d}\:\sum_\pm\sum_{\substack{m',m''|cm\\m'\sim \mathfrak{m}'}}\:\sum_{n=1}^\infty \frac{\lambda(n,m'')}{m''n}\\
\nonumber &\times\sideset{}{^\star}\sum_{\substack{\beta\bmod{mc/m'}\\r+\beta m'\equiv 0\bmod{c/d}}}S(m'\beta, \pm n; mc/m'')\;V^\star_\pm\left(\frac{c}{C},\frac{m'^2m''^2nL}{m^3c^3},\frac{r}{\tilde{R}}\right).
\end{align}
The integral transform is negligibly small if 
\begin{align}
\label{mathfrakn}
n>\mathfrak{N}=\frac{m^3C^3M^{3\hat\theta+\varepsilon}}{(m'm'')^2L}.
\end{align}
Substituting the above expression (only the term with $+$ sign) in the second line of \eqref{dual-od-red-tran-2} we get
\begin{align*}
&\frac{RN\tilde{N}}{m^2CMP^5}\sum_{\substack{P<p<2P\\p\;\text{prime}}}\;\phi(p)\chi(p)
\mathop{\sum}_{c=1}^\infty\sum_{\substack{r=1\\(p,r)=1}}^\infty \bar{\chi}(r)e\left(\frac{Mc\overline{r}}{p}\right)\sum_{\substack{d|c\\d\sim \mathfrak{d}}}\frac{\mu(d)}{d}\:\sum_{\substack{m',m''|cm\\m'\sim\mathfrak{m}'}}\\
&\sum_{n\leq\mathfrak{N}} \frac{\lambda(n,m'')}{m''n}\sideset{}{^\star}\sum_{\substack{\beta\bmod{mc/m'}\\r+\beta m'\equiv 0\bmod{c/d}}}S(m'\beta, n; mc/m'')\;V_+^\star\left(\frac{c}{C},\frac{m'^2m''^2nL}{m^3c^3},\frac{r}{\tilde{R}}\right).
\end{align*}
Here the sum over $n$ is truncated at $\mathfrak{N}$ at a cost of a negligible error term. Taking inverse Mellin transform to free the variable $r$ from the weight function, and then taking absolute values we bound the above sum by
\begin{align*}
&\frac{RN\tilde{N}}{m^2CMP^5}
\int_{-M^{\hat{\theta}+\varepsilon}}^{M^{\hat{\theta}+\varepsilon}}\mathop{\sum}_{C<c<2C}\:\sum_{\substack{d|c\\d\sim\mathfrak{d}}}\frac{1}{d}\:\sum_{\substack{m',m''|cm\\m'\sim\mathfrak{m}'}}\;\sum_{\tilde{R}<r<2\tilde{R}} \Bigl|\sum_{\substack{P<p<2P\\p\;\text{prime}\\(p,r)=1}}\;\phi(p)\chi(p)e\left(\frac{Mc\overline{r}}{p}\right)\Bigr|\\
&\times \Bigl|\sum_{n\leq\mathfrak{N}} \frac{\lambda(n,m'')}{m''n}\sideset{}{^\star}\sum_{\substack{\beta\bmod{mc/m'}\\r+\beta m'\equiv 0\bmod{c/d}}}S(m'\beta, n; mc/m'')\tilde{V}_+^\star\left(\frac{c}{C},\frac{m'^2m''^2nL}{m^3c^3},iu\right)\Bigr|\mathrm{d}u.
\end{align*}
Recall that the weight function $V_+^\star$ does not depend on $p$.\\

Applying Cauchy inequality we get that the above sum is dominated by 
\begin{align*}
M^{\hat{\theta}+\varepsilon}\frac{RN\tilde{N}}{m^2CMP^5}\mathfrak{d}^{-1/2}\:\sqrt{\Theta_1}\:\sqrt{\Theta_3(\mathfrak{m}',\mathfrak{d})}
\end{align*}
where $\Theta_1$ is as given in \eqref{theta1}, and 
\begin{align*}
&\Theta_3(\mathfrak{m}',\mathfrak{d})=\sup_u \mathop{\sum}_{C<c<2C}\:\sum_{\substack{d|c\\d\sim \mathfrak{d}}}\frac{1}{d}\:\sum_{\substack{m',m''|cm\\m'\sim\mathfrak{m}'}}\;\sum_{r\in\mathbb{Z}}W\left(\frac{r}{\tilde{R}}\right)\\
\times \Bigl|\sum_{n\leq \mathfrak{N}} &\frac{\lambda(n,m'')}{m''n}\sideset{}{^\star}\sum_{\substack{\beta\bmod{mc/m'}\\r+\beta m'\equiv 0\bmod{c/d}}}S(m'\beta, n; mc/m'')\tilde V_+^\star\left(\frac{c}{C},\frac{m'^2m''^2nL}{m^3c^3},iu\right)\Bigr|^2.
\end{align*}
Here the supremum is taken over the range $|u|\ll M^{\hat{\theta}+\varepsilon}$ and $W$ is a non-negative compactly supported smooth function on $(0,\infty)$ with $W(x)=1$ for $x\in [1,2]$. We conclude, using \eqref{bd-theta1}, that
\begin{align}
\label{mid-bd-dual}
\mathcal{O}^+(C,m;L,\mathfrak{m}',\mathfrak{d})\ll M^{\hat{\theta}+\varepsilon}\frac{RN\tilde{N}}{m^2CMP^5}\:\mathfrak{d}^{-1/2}\:\sqrt{P^3C\tilde{R}}\:\sqrt{\Theta_3(\mathfrak{m}',\mathfrak{d})}.
\end{align}
Here the $+$ denotes that we are only considering the contribution of the $+$ term from \eqref{again-voronoi}. The contribution of the $-$ term can be analysed in the similar fashion. The bound that we obtain is not sensitive to this sign. \\

\subsection{Bound for $\Theta_3$}
Now we consider $\Theta_3=\Theta_3(\mathfrak{m}',\mathfrak{d})$. 
Opening the absolute square we perform Poisson summation on the $r$ sum with modulus $c/d$. We get
\begin{align}
\label{theta3-poi}
\Theta_3=\sup_u\:\tilde{R}\mathop{\sum}_{C<c<2C}\:\frac{1}{c}\sum_{\substack{d|c\\d\sim\mathfrak{d}}}&\:\sum_{\substack{m',m''|cm\\m'\sim\mathfrak{m}'}}\:\mathop{\sum\sum}_{n,n'\leq\mathfrak{N}} \frac{\lambda(n,m'')}{m''n}\frac{\overline{\lambda(n',m'')}}{m''n'}\tilde{V}_+^\star(\dots)\bar{\tilde{V}}_+^\star(\dots) \sum_{r\in\mathbb{Z}}\:\mathfrak{I}\:\mathfrak{C}
\end{align}
where the character sum is given by
\begin{align*}
\mathfrak{C}=\mathop{\sideset{}{^\star}\sum\sideset{}{^\star}\sum}_{\substack{\beta,\beta'\bmod{mc/m'}\\\beta m'\equiv \beta'm'\bmod{c/d}}}S(m'\beta, n; mc/m'')S(m'\beta', n'; mc/m'')e\left(-\frac{r\beta m'}{c/d}\right)
\end{align*}
and the integral is given by
\begin{align*}
\mathfrak{I}=\int_{\mathbb{R}}W(z)e\left(-\frac{z\tilde{R}r}{c/d}\right)\mathrm{d}z.
\end{align*}
By repeated integration by parts it follows that the integral is negligibly small if 
$$
|r|\gg \frac{CM^{\varepsilon}}{\mathfrak{d}\tilde{R}}.
$$
Given the restriction on the sizes of $C$ and $\tilde{R}$, we see that the non-zero frequencies $r\neq 0$ make a negligible contribution. For $r=0$ we use the trivial bound for the integral $\mathfrak{I}\ll 1$.\\ 

Next we will estimate the character sum for $r=0$. We temporarily denote this by $\mathfrak{C}_0$. Let $p$ be a prime with $v_p(c)=\gamma$, $v_p(d)=\delta$, $v_p(m)=\mu$, $v_p(m')=\mu'$ and $v_p(m'')=\mu''$. We consider the character sum
\begin{align*}
\mathop{\sideset{}{^\star}\sum\sideset{}{^\star}\sum}_{\substack{\beta,\beta'\bmod{p^{\gamma+\mu-\mu'}}\\\beta p^{\mu'}\equiv \beta'p^{\mu'}\bmod{p^{\gamma-\delta}}}}S(p^{\mu'}\beta a, nb; p^{\gamma+\mu-\mu''})S(p^{\mu'}\beta' a, n' b; p^{\gamma+\mu-\mu''})
\end{align*}
where $p\nmid ab$. If $\mu'\geq\gamma-\delta$, then the sum splits into a product of two sums
\begin{align*}
\mathop{\sideset{}{^\star}\sum}_{\substack{\beta\bmod{p^{\gamma+\mu-\mu'}}}}S(p^{\mu'}\beta a, nb; p^{\gamma+\mu-\mu''})\mathop{\sideset{}{^\star}\sum}_{\substack{\beta'\bmod{p^{\gamma+\mu-\mu'}}}}S(p^{\mu'}\beta' a, n' b; p^{\gamma+\mu-\mu''}),
\end{align*}
which can be written as a product of Ramanujan sums
$$
\mathfrak{c}_{p^{\gamma+\mu-\mu''}}(n)\mathfrak{c}_{p^{\gamma+\mu-\mu''}}(n')\mathfrak{c}^2_{p^{\gamma+\mu-\mu'}}(p^{\mu''})\leq (p^{\gamma+\mu-\mu''},n)(p^{\gamma+\mu-\mu''},n')(p^{\gamma+\mu-\mu'},p^{\mu''})^2.
$$ 
The last term can be bounded by
$$
p^{\gamma+\mu-\mu'+\mu''}\Bigl((p^{\gamma+\mu-\mu''},n)+(p^{\gamma+\mu-\mu''},n')\Bigr)(p^{\gamma+\mu-\mu''},n-n').
$$\\

On the other hand if $\mu'<\gamma-\delta$, then we have congruence restriction $\beta\equiv\beta'\bmod{p^{\gamma-\delta-\mu'}}$, and the above sum boils down to
\begin{align*}
\mathop{\sideset{}{^\star}\sum}_{\beta\bmod{p^{\gamma-\delta-\mu'}}}\sum_{\beta_1\bmod{p^{\mu+\delta}}}&S\Bigl(p^{\mu'}(\beta+\beta_1p^{\gamma-\delta-\mu'}) a, nb; p^{\gamma+\mu-\mu''}\Bigr)\\
&\sum_{\beta_1'\bmod{p^{\mu+\delta}}}S\Bigl(p^{\mu'}(\beta+\beta_1'p^{\gamma-\delta-\mu'}) a, n' b; p^{\gamma+\mu-\mu''}\Bigr).
\end{align*}
Opening the Kloosterman sums we observe that the sums over $\beta_1$ and $\beta_1'$ vanishes unless $\mu''=\mu+\delta$. In this case we also need $p^{\mu'}|n$ and $n'$, otherwise the average of the Kloosterman sum vanishes. Set $n=p^{\mu'}\tilde{n}$ and $n'=p^{\mu'}\tilde{n}'$. The character sum now reduces to
\begin{align*}
p^{2(\mu+\delta)}\mathop{\sideset{}{^\star}\sum}_{\beta\bmod{p^{\gamma-\delta-\mu'}}}S\Bigl(\beta a, \tilde{n}b; p^{\gamma-\delta-\mu'}\Bigr)S\Bigl(\beta a, \tilde{n}' b; p^{\gamma-\delta-\mu'}\Bigr).
\end{align*}
As $p\nmid ab$ we can change variables to arrive at
\begin{align*}
p^{2(\mu+\delta)}\mathop{\sideset{}{^\star}\sum}_{\beta\bmod{p^{\gamma-\delta-\mu'}}}S\Bigl(\beta , \tilde{n}; p^{\gamma-\delta-\mu'}\Bigr)S\Bigl(\beta , \tilde{n}' ; p^{\gamma-\delta-\mu'}\Bigr)
\end{align*}
which is given by
\begin{align*}
p^{\gamma+2\mu+\delta-\mu'}\mathfrak{c}_{p^{\gamma-\delta-\mu'}}(\tilde{n}-\tilde{n}')-p^{\gamma+2\mu+\delta-\mu'-1}\mathfrak{c}_{p^{\gamma-\delta-\mu'-1}}\left(\frac{\tilde{n}}{p}-\frac{\tilde{n}'}{p}\right).
\end{align*}
This is bounded by
$$
p^{\gamma+\mu-\mu'+\mu''}(p^{\gamma+\mu-\mu''},n-n').
$$\\

\begin{lemma}
We have
\begin{align}
\label{mathfrakC}
\mathfrak{C}_0\ll \frac{cmm''}{m'}\Bigl((cm/m'',n)+(cm/m'',n')\Bigr)(cm/m'',n-n').
\end{align}
\end{lemma}

\bigskip

It now follows that
\begin{align*}
\Theta_3\ll \sup_u\:\tilde{R}\mathop{\sum}_{C<c<2C}\:\frac{1}{c}\sum_{\substack{d|c\\d\sim\mathfrak{d}}}&\:\sum_{\substack{m',m''|cm\\m'\sim\mathfrak{m}'}}\:\mathop{\sum\sum}_{n,n'\leq\mathfrak{N}} \frac{|\lambda(n,m'')|}{m''n}\frac{|\lambda(n',m'')|}{m''n'}|\tilde{V}_+^\star(\dots)||\tilde{V}_+^\star(\dots)|\\
&\times \frac{cmm''}{m'}(cm/m'',n)(cm/m'',n-n')+M^{-2013}.
\end{align*}
Here it is not clear whether one can estimate this sum without taking point wise bound for the Fourier coefficients. Direct application of Cauchy is not helpful as the gcd function has a large dispersion. Using (R) and (RS) we get
\begin{align}
\label{mathfrakI}
\frac{\lambda(n,m'')}{m''n}\frac{\overline{\lambda(n',m'')}}{m''n'}\tilde{V}_+^\star(\dots)\bar{\tilde{V}}_+^\star(\dots)\ll \frac{1}{nn'm''^2}\frac{m'^2m''^2nL}{m^3c^3}\frac{m'^2m''^2n'L}{m^3c^3}M^\varepsilon.
\end{align}
Substituting this in the above expression, we get
\begin{align*}
\Theta_3\ll M^\varepsilon\tilde{R}\mathop{\sum}_{C<c<2C}\:\sum_{\substack{d|c}}&\:\sum_{\substack{m',m''|cm\\m'\sim\mathfrak{m}'}}\:\frac{m'^4m''^3L^2}{c^6m^6}\frac{m}{m'}\:\mathop{\sum\sum}_{1\leq n,n'\leq\mathfrak{N}}\:(cm/m'',n)(c/m'',n-n').
\end{align*}\\

Next we sum over $n$ and $n'$. The contribution from the diagonal $n=n'$ is dominated by $c\mathfrak{N}/m''$ and the off-diagonal is dominated by $\mathfrak{N}^2$. Hence 
\begin{align*}
\Theta_3\ll \frac{mM^\varepsilon\tilde{R}L^2}{\mathfrak{m}'}\mathop{\sum}_{C<c<2C}\:\sum_{d|c}&\:\sum_{\substack{m',m''|cm\\m'\sim\mathfrak{m}'}}\:\frac{m'^4m''^3}{(cm)^6}\left(\frac{c}{m''}\mathfrak{N}+\mathfrak{N}^2\right).
\end{align*}
Substituting the size of $\mathfrak{N}$ from \eqref{mathfrakn}, we get 
\begin{align*}
\Theta_3\ll \frac{mM^{\varepsilon}\tilde{R}L^2}{\mathfrak{m}'}\mathop{\sum}_{C<c<2C}\:\sum_{d|c}&\:\sum_{\substack{m',m''|cm\\m'\sim\mathfrak{m}'}}\:\frac{m'^4m''^3}{(cm)^6}\left(\frac{m^3C^4M^{3\hat\theta}}{m'^2m''^3L}+\frac{m^6C^6M^{6\hat\theta}}{(m'm'')^4L^2}\right).
\end{align*}
Now applying the upper bound for $L$ from \eqref{l-range} we arrive at
\begin{align*}
\Theta_3\ll &\frac{mM^{\varepsilon}\tilde{R}}{\mathfrak{m}'}\mathop{\sum}_{C<c<2C}\:\sum_{\substack{m',m''|cm\\m'\sim\mathfrak{m}'}}\:\frac{m'^4m''^3}{(cm)^6}\\
&\times \left(\frac{m^3C^4M^{3\hat\theta}}{m'^2m''^3}\frac{C^3m^3M^{3\hat \theta}}{\mathfrak{m}'^2\tilde{N}}+\frac{m^6C^6M^{6\hat\theta}}{(m'm'')^4}\right).
\end{align*}

\begin{lemma}
We have
\begin{align*}
\Theta_3\ll \frac{mM^{6\hat\theta + \varepsilon}\tilde{R}}{\mathfrak{m}'}\left(\frac{C^2}{\tilde{N}}+C\right).
\end{align*}
\end{lemma}

\bigskip

\subsection{Conclusion}
Substituting the above bound for $\Theta_3$ in \eqref{mid-bd-dual} we get
\begin{align*}
&\mathcal{O}^+(C,m;L,\mathfrak{m}',\mathfrak{d})
\ll  M^\varepsilon\frac{RN\tilde{N}}{m^2CMP^5}\:\mathfrak{d}^{-1/2}\:\sqrt{P^3C\tilde{R}}\:\sqrt{\frac{mM^{6\hat\theta}\tilde{R}}{\mathfrak{m}'}\left(\frac{C^2}{\tilde{N}}+C\right)}\\
\ll &\frac{M^{3\hat\theta+\varepsilon}}{m^{3/2}\sqrt{\mathfrak{m}'\mathfrak{d}}}\frac{RN\sqrt{\tilde{R}\tilde{N}}}{MP^5}\sqrt{P^3\tilde{R}(C+\tilde{N})}\ll \frac{M^{3\hat\theta+\varepsilon}}{m^{3/2}\sqrt{\mathfrak{m}'\mathfrak{d}}}M^{3/2}\left(\frac{N}{R}\right)^{1/4}.
\end{align*}
The following lemma summarizes the main content of this section.\\

\begin{lemma}
\label{lem-62}
Suppose $N^{1/2+\varepsilon}<P<M^{1-\varepsilon}$ and $\mathfrak{m}'\mathfrak{d}\geq M^{28\theta}$. Then we have 
\begin{align*}
\mathcal{O}(C,m;L,\mathfrak{m}',\mathfrak{d})\ll m^{-3/2}\sqrt{N}M^{3/4-\theta/2+\varepsilon}.
\end{align*}\\
\end{lemma}

Combining with Lemma~\ref{lemma15} we draw the following conclusion.\\

\begin{corollary}
\label{cor-lemma15}
We have
\begin{align*}
\mathcal{O}^\star(C,m)\ll M^\varepsilon \sup |\mathcal{O}(C,m;L,\mathfrak{m}',\mathfrak{d})|+m^{-3/2}\sqrt{N}M^{3/4-\theta/2+\varepsilon}
\end{align*}
where the supremum is taken over all $L$ in the range \eqref{l-range}, $\mathfrak{m}'\mathfrak{d}\ll M^{28\theta}$.\\
\end{corollary}

Combining with Corollary~\ref{cor-previous-11} and Lemma~\ref{lemma-nine}, we conclude the following.\\
 
\begin{corollary}
\label{cor-previous-1115}
Suppose $N^{1/2+\varepsilon}<P<M^{1-\varepsilon}$, and $\theta<1/8$, then we have
$$
\mathcal{F}\ll M^\varepsilon\sum_{m\leq M^{4\theta}}\:\sup\:|\mathcal{O}(C,m;L,\mathfrak{m}',\mathfrak{d})|+\sqrt{N}M^{3/4-\theta/2+\varepsilon}
$$
where the supremum is taken over all $L$ in the range \eqref{l-range}, $\mathfrak{m}'\mathfrak{d}\ll M^{28\theta}$, and $C$ in the range \eqref{range-for-C}.
\end{corollary}

\bigskip

\section{Tamed dual off-diagonal in transition}
\label{dod3}

We now return to \eqref{dual-od-red-tran-2} the expression we obtained after the first application of the Voronoi summation and dyadic segmentation. Let $\theta^\star=28\theta$ and $\hat\theta=4\theta$. We take $C$ in the transition range \eqref{range-for-C}, $m$ in the range $1\leq m\leq M^{\hat\theta}$. We write $cd$ in place of $c$ and change the order of summations. It follows that 
\begin{align}
\label{bd-inter}
\sup_{\mathfrak{m}'\mathfrak{d}\leq M^{\theta^\star}}\mathcal{O}(C,m;L,\mathfrak{m}',\mathfrak{d})\ll M^\varepsilon\sup_{d\mathfrak{m}'\leq M^{\theta^\star}}\Bigl|\mathfrak{O}(\dots)\Bigr|
\end{align} 
where
\begin{align}
\label{dual-od-red-tran-3}
&\mathfrak{O}(\dots)=\frac{RN\tilde{N}_0}{CMP^5}\sum_{\substack{P<p<2P\\p\;\text{prime}}}\;\phi(p)\chi(p)\sum_{\substack{r=1\\(p,r)=1}}^\infty \bar{\chi}(r)\mathop{\sum}_{c=1}^\infty e\left(\frac{Mcd\overline{r}}{p}\right)\\
\nonumber &\times \sum_{\substack{m'|cdm\\m'\sim \mathfrak{m}'}}\frac{m'}{cdm}\:\sideset{}{^\star}\sum_{\substack{\beta\bmod{mcd/m'}\\r+\beta m'\equiv 0\bmod{c}}}\;\sum_{n=1}^\infty \lambda(m',n)e\left( \frac{\bar{\beta}n}{mcd/m'}\right)V\left(\frac{cd}{C},\frac{n}{L},\frac{r}{\tilde{R}}\right).
\end{align}
Here the weight function $V$ is as given in \eqref{after-voronoi}.\\

\subsection{Evaluation of character sum and reciprocity}
Consider the character sum (which we again temporarily denote by $\mathfrak{C}$)
\begin{align*}
\mathfrak{C}=\sideset{}{^\star}\sum_{\substack{\beta\bmod{mcd/m'}\\r+\beta m'\equiv 0\bmod{c}}}e\left( \frac{\bar{\beta}n}{mcd/m'}\right).
\end{align*}
If $m=m'=d=1$, then the character sum can be explicitly evaluated and it is given by $e(-\bar{r}n/c)$. However in general it is not easy to evaluate the character sum due to the presence of factors $m$, $m'$ and $d$. But we have now obtained a good control on the sizes of these factors, and consequently we can evaluate explicitly a large `portion' of the character sum.\\
 
To this end let $h=(m',c)$. We observe that $\mathfrak{C}=0$ unless $h|r$. Accordingly we write $m'=hm_1'$, $c=hc_1$ and $r=hr_1$. Let $h_1=(m_1',r_1)$ and let us write $r_1=h_1r_2$ and $m_1'=h_1m_2$. Hence $(r_2,m_2)=1$. We get
\begin{align*}
\mathfrak{C}=\sideset{}{^\star}\sum_{\substack{\beta\bmod{mc_1d/m_1'}\\ \beta\equiv -r_2\overline{m}_2\bmod{c_1}}}e\left( \frac{\bar{\beta}n}{mc_1d/m_1'}\right).
\end{align*}
It follows that
\begin{align*}
\mathfrak{O}(\dots)=&\frac{RN\tilde{N}_0}{CMP^5}\mathop{\sum\sum}_{\substack{h_1,m_2\\h_1m_2=m_1'|dm}}\frac{m_1'}{dm}\;\sum_{h\sim \mathfrak{m}'/m_1'}\bar{\chi}(hh_1)\\
\times &\sum_{\substack{P<p<2P\\p\;\text{prime}}}\;\phi(p)\chi(p)
\sum_{\substack{r_2=1\\(pm_2,r_2)=1}}^\infty \bar{\chi}(r_2)\mathop{\sum}_{\substack{c_1=1\\(r_2m_1',c_1)=1}}^\infty e\left(\frac{Mc_1d\overline{h_1r_2}}{p}\right)\frac{1}{c_1}\\
\times &\sum_{n=1}^\infty \lambda(m',n)\sideset{}{^\star}\sum_{\substack{\beta\bmod{mc_1d/m_1'}\\ \beta\equiv -r_2\overline{m}_2\bmod{c_1}}}e\left( \frac{\bar{\beta}n}{mc_1d/m_1'}\right)V\left(\frac{c_1hd}{C},\frac{n}{L},\frac{r_2hh_1}{\tilde{R}}\right)
\end{align*}
where $m'=hm_1'$.\\ 

Let $g=(c_1,dm)$. We write $c_1=gc_2$ and $dm=g_0g'$, where $g_0|g^\infty$ and $(g',g)=1$. Let $f=g'/m_1'$, which is an integer as $m_1'|g_0g'$ but $(c_1,m_1')=1$. Then $\bar\beta=-\bar{r}_2m_2+\beta_1 c_1$ with $\beta_1 \bmod{fg_0}$. We have
$$
\mathfrak{C}=e\left( -\frac{\bar{r}_2m_2n}{c_1fg_0}\right)\sideset{}{^\dagger}\sum_{\substack{\beta_1\bmod{fg_0}}}e\left( \frac{\beta_1n}{fg_0}\right)
$$
where $\dagger$ implies that $(\beta_1,f)=1$. In particular $\mathfrak{C}=0$ if $(r_2,c_1fg_0)>1$. Applying the reciprocity relation to the outer exponential, and pulling out the gcd of $\beta_1$ and $g_0$, we get
$$
\mathfrak{C}=e\left(\frac{\overline{c_1fg_0}m_2n}{r_2}\right)e\left( -\frac{m_2n}{c_1fg_0r_2}\right)\sum_{g_1g_2=g_0}g_1\sum_{\substack{1\leq \beta_1<fg_2\\(\beta_1,fg_2)=1}}e\left( \frac{\beta_1n}{fg_2}\right).
$$
Rearranging sums we arrive at
\begin{align}
\label{ohhm}
\mathfrak{O}(\dots)\ll &\frac{RN\tilde{N}_0}{CMP^5}\mathop{\sum\sum}_{\substack{h_1,m_2\\h_1m_2=m_1'|dm}}\frac{m_1'}{dm}\sum_{\substack{g_0g'=dm\\(g_0,g')=1\\(g_0,m_1')=1}}\sum_{g|g_0|g^\infty}\\
\nonumber &\times \sum_{h\sim \mathfrak{m}'/m_1'}\sum_{g_1g_2=g_0}\frac{g_1}{g}\sum_{\substack{1\leq \beta_1<fg_2\\(\beta_1,fg_2)=1}}|\Omega(\dots)|
\end{align}
where $f=dm/g_0m_1'$, 
\begin{align*}
\Omega(\dots)=\sum_{\substack{P<p<2P\\p\;\text{prime}}}&\;\phi(p)\chi(p)\sum_{\substack{r_2=1\\(fgpm_2,r_2)=1}}^\infty \bar{\chi}(r_2)\mathop{\sum}_{\substack{c_2=1\\(fg_3m_1'r_2,c_2)=1}}^\infty e\left(\frac{Mgc_2d\overline{h_1r_2}}{p}\right)\frac{1}{c_2}\\
&\sum_{n=1}^\infty \lambda(m',n)e\left(\frac{\overline{c_2fgg_0}m_2n}{r_2}+ \frac{\beta_1n}{fg_2}\right)W\left(\frac{c_2ghd}{C},\frac{n}{L},\frac{r_2hh_1}{\tilde{R}}\right)
\end{align*}
with $g_3=g_0/g$ and
\begin{align*}
W(x,y,z)=e\left( -\frac{m_2h^2h_1d}{fg_0}\frac{L}{C\tilde{R}}\frac{y}{xz}\right)V(x,y,z).
\end{align*}
 The new weight function $W$ is smooth, supported in $[1,2]^3$ and satisfies 
\begin{align}
\label{derivatives}
W^{(j_1,j_2,j_3)}(x,y,z)\ll M^{12\theta(j_1+j_2+j_3)}. 
\end{align} 
\\

\begin{lemma}
\label{lem-B}
Suppose
$$
B(\dots)=\sup\: |\Omega(\dots)|
$$
where the supremum is taken over all possible vectors $(d,m,h_1,m_2,g,h,g_1,g_2,\beta_1)$ which appear in \eqref{ohhm}, with $m\leq M^{\hat\theta}$, and $dhh_1m_2\leq M^{\theta^\star}$. Then
\begin{align*}
\sup_{\mathfrak{m}'\mathfrak{d}\leq M^{\theta^\star}}\:|\mathcal{O}(C,m;L,\mathfrak{m}',\mathfrak{d})| \ll M^{\theta^\star+\varepsilon}\:B(\dots)\:\frac{RN\tilde{N}_0}{CMP^5}.
\end{align*} 
\end{lemma}

\bigskip

\subsection{The last application of Voronoi summation}
In the rest of the paper we will obtain a sufficient bound for $B(\dots)$. We apply the Voronoi summation on the sum over $n$ in $\Omega(\dots)$. The modulus of the additive character is $r_2fg_2$. Notice that the application of the reciprocity relation has changed the modulus and so the Voronoi summation here is not a reversal process. This gives rise to two terms - a $+$ term and a $-$ term. We will analyse the contribution of the $+$ term, which is given by
\begin{align}
\label{omega+}
&\Omega_+(\dots)=\sum_{\substack{P<p<2P\\p\;\text{prime}}}\;\phi(p)\chi(p)\sum_{\substack{r_2=1\\(fgpm_2,r_2)=1}}^\infty \bar{\chi}(r_2)\\
\nonumber &\times \mathop{\sum}_{\substack{c_2=1\\(fg_3m_1'r_2,c_2)=1}}^\infty e\left(\frac{Mgc_2d\overline{h_1r_2}}{p}\right)\frac{fg_2r_2}{c_2}\:\sum_{m''|fg_2r_2m'}\\
\nonumber &\times \sum_{n=1}^\infty \frac{\lambda(m'',n)}{m''n}S(m'\bar{\xi}, n; m'fg_2r_2/m'')W_+^\star\left(\frac{c_2ghd}{C},\frac{m''^2nL}{(fg_2r_2)^3m'},\frac{r_2hh_1}{\tilde{R}}\right)
\end{align}
where $\xi=\overline{c_2fgg_0}fg_2m_2+\beta_1r_2$. (Note that $\xi$ is invertible modulo $fg_2r_2$, as $(fg_2,r_2)=1$, and $(\beta_1,fg_2)=(c_2fgm_2,r_2)=1$.) The integral transform is negligibly small if
\begin{align}
\label{new-N}
n\geq \mathcal{N}=\frac{(fg_2\tilde{R})^3m'M^{36\theta+\varepsilon}}{(hh_1)^3m''^2L}.
\end{align}
For smaller values of $n$ we shift the contour in the definition of the integral transform \eqref{gl} to $\sigma=\varepsilon$, using (RS). The integrand decays rapidly for $t=\text{Im}(s)\gg M^{12\theta+\varepsilon}$ and this part makes a negligible contribution. We now interchange the order of summations and the integral over $t$. This reduces the analysis of the above sum $\Omega_+$ to sums of the form
\begin{align}
\label{41}
L&\sum_{\substack{P<p<2P\\p\;\text{prime}}}\;\phi(p)\chi(p)\sum_{\substack{r_2=1\\(fgpm_2,r_2)=1}}^\infty \bar{\chi}(r_2)\\
\nonumber &\times \mathop{\sum}_{\substack{c_2=1\\(fg_3m_1'r_2,c_2)=1}}^\infty e\left(\frac{Mgc_2d\overline{h_1r_2}}{p}\right)\frac{1}{c_2}\sum_{m''|fg_2r_2m'}\:\frac{m''}{(fg_2r_2)^2m'}\\
\nonumber &\times \sum_{1\leq n<\mathcal{N}} \frac{\lambda(m'',n)}{n^{it}}\:S(m'\bar{\xi}, n; m'fg_2r_2/m'')\: U_t\left(\frac{c_2ghd}{C},\frac{r_2hh_1}{\tilde{R}}\right).
\end{align}
The weight function $U_t$ is smooth, supported in $[1,2]^2$ and satisfies
$$
U_t^{(i,j)}(x,y)\ll M^{12\theta(i+j)},
$$
where the implied constant is independent of $t$. In the rest of the paper we will obtain sufficient bounds for the expression in \eqref{41}, which is uniform with respect $t$ in the desired range. Such a bound when multiplied by $M^{12\theta+\varepsilon}$ will yield a bound for $\Omega(\dots)$.\\

\subsection{Reciprocity and Poisson summation}
Next we wish to apply the Poisson summation formula on the sum over $c_2$. Recall that $c_2$ is essentially the variable $c$ which is the modulus of the `circle method' (Petersson formula) that we applied at the initial stage. After a sequence of applications of summation formulas and reciprocity relations we are finally at the stage where we are able to sum over the modulus again. The variable $c_2$ appears in the Kloosterman sum in \eqref{omega+}. This Kloosterman sum has modulus $fg_2m'r_2/m''$. Also $c_2$ appears in the additive character which has modulus $p$. So apparently the total modulus is too large compared to the length of the sum. However we can now apply the reciprocity relation 
$$
e\left(\frac{Mgc_2d\overline{h_1r_2}}{p}\right)=e\left(-\frac{Mgc_2d\overline{p}}{h_1r_2}\right)e\left(\frac{Mgc_2d}{h_1pr_2}\right).
$$
The last term can be absorbed in the weight function. Accordingly we let
$$
V_t(x,y)=e\left(\frac{CM}{p\tilde{R}}\frac{x}{y}\right)U_t(x,y).
$$
Observe that we (still) have
$$
V_t^{(i,j)}(x,y)\ll M^{12\theta(i+j)}.
$$\\

We now study the sum over $c_2$ in \eqref{41} which is given by
\begin{align*}
\mathop{\sum}_{\substack{c_2=1\\(fg_3m_1'r_2,c_2)=1}}^\infty e\left(-\frac{Mgc_2d\overline{p}}{h_1r_2}\right)\frac{1}{c_2}\:S(m'\bar{\xi}, n; m'fg_2r_2/m'')\: V_t\left(\frac{c_2ghd}{C},\frac{r_2hh_1}{\tilde{R}}\right).
\end{align*}
We break the sum into congruence classes modulo $fg_0m'r_2=fg_0hm_1'r_2=dmhr_2$ and apply the Poisson summation formula. We get
\begin{align*}
\frac{1}{dmhr_2}\sum_{c_2\in\mathbb{Z}}\mathop{\sum}_{\substack{\gamma\bmod{dmhr_2}\\(fg_3m_1'r_2,\gamma)=1}} \:&e\left(-\frac{Mgd\overline{p}\gamma}{h_1r_2}+\frac{c_2\gamma}{dmhr_2}\right)\:S(m'\bar{\xi}, n; m'fg_2r_2/m'')\\
&\times \int_{\mathbb{R}}V_t\left(x,\frac{r_2hh_1}{\tilde{R}}\right)e\left(-\frac{Cc_2x}{d^2gh^2mr_2}\right)\frac{\mathrm{d}x}{x},
\end{align*}
where $\xi=\overline{\gamma fgg_0}fg_2m_2+\beta_1r_2$. From repeated integration by parts it follows that the integral is negligibly small if 
$$
|c_2|\gg \mathcal{C}_2=M^{12\theta+\varepsilon}\frac{d^2ghm\tilde{R}}{Ch_1}.
$$\\

\subsection{Evaluation of character sums}
Now we write $r_2=r_3r_4$ with $(r_3,dmh)=1$ and $r_4|(dmh)^\infty$. Accordingly we split $m''=m_3''m_4''$, with $m_3''|r_3$ and $m_4''|fg_2m'r_4$. The character sum now splits as a product of two character sums. The one with modulus $dmhr_4$ is given by
\begin{align*}
\sum_{\substack{\gamma\bmod{dmhr_4}\\(fg_3m_1'r_4,\gamma)=1}}S(m'\overline{\xi r_3/m_3''}, n\overline{r_3/m_3''}; m'fg_2r_4/m_4'')e\left(-\frac{Mg\gamma d\overline{pr_3}}{h_1r_4}+\frac{\gamma c_2\overline{r_3}}{dmhr_4}\right).
\end{align*}
Suppose $p^\ell\|r_4$ with $\ell\geq 1$, and suppose $p^k\|dmh$, $p^j\|m'fg_2r_4/m_4''$ (so $j\leq \ell+k$). Then we study the sum
\begin{align*}
\sideset{}{^\star}\sum_{\gamma\bmod{p^{\ell+k}}}S(m'A\overline{\xi}, B; p^j)e\left(\frac{C\gamma}{p^{\ell+k}}\right)
\end{align*}
where $p\nmid A$. The sum vanishes unless $p^{\ell+k-j}|C$, in which case it reduces to
\begin{align*}
p^{\ell+k-j}\sideset{}{^\star}\sum_{\gamma\bmod{p^j}}S(m'A\overline{\xi}, B; p^j)e\left(\frac{C\gamma}{p^j}\right).
\end{align*}
First consider the case where $\ell> k$, so that $2\ell> j$. Then 
$$
\bar{\xi}=\gamma fgg_0\overline{fg_2m_2}-(\gamma fgg_0\overline{fg_2m_2})^2\beta_1r_2.
$$
Now if $\ell\geq j$ the character sum reduces to
\begin{align*}
p^{\ell+k-j}\sideset{}{^\star}\sum_{\gamma\bmod{p^j}}S(m'A\gamma gg_0\overline{g_2m_2}, B; p^j)e\left(\frac{C\gamma}{p^j}\right).
\end{align*}
Opening the Kloosterman sum we execute the sum over $\gamma$, which yields a Ramanujan sum. Using standard bounds for the Ramanujan sum we now get the bound $O(p^{\ell+k}(m',p^j))$ for the character sum. On the other hand if $j>\ell\geq k$, then we write $\gamma=\gamma_1+\gamma_2p^{j-\ell}$ with $\gamma_1$ modulo $p^{j-\ell}$, $(\gamma_1,p)=1$, and $\gamma_2$ modulo $p^\ell$. Then the character sum reduces to
\begin{align*}
p^{\ell+k-j}\mathop{\sideset{}{^\star}\sum\sum}_{\substack{\gamma_1\bmod{p^{j-\ell}}\\\gamma_2\bmod{p^\ell}}}&S(m'A((\gamma_1 +\gamma_2p^{j-\ell}gg_0\overline{g_2m_2})-(\gamma_1 gg_0\overline{g_2m_2})^2\beta_1r_2), B; p^j)\\
&\times e\left(\frac{C\gamma_1}{p^j}+\frac{C\gamma_2}{p^\ell}\right).
\end{align*}
Opening the Kloosterman sum, executing the sum over $\gamma_2$, and trivially estimating the remaining sums we get the bound $O(p^{j+k}(m',p^\ell))$. In the case $\ell<k$ (including when $\ell=0$) we trivially bound the sum by $O(p^{\ell+k+j})=O(p^{2k+j})$. Putting the above bounds together we are able to bound the initial character sum with modulus $dmhr_4$ by
$$
O\left((dmh)^2\frac{fg_2m'r_4}{m_4''}\right).
$$\\

The other character sum with modulus $r_3$ is given by
\begin{align*}
\sideset{}{^\star}\sum_{\gamma\bmod{r_3}}S(m'\overline{\eta}, n; r_3/m_3'')e\left(-\frac{Mg\gamma d\overline{ph_1r_4}}{r_3}+\frac{\gamma c_2\overline{dmhr_4}}{r_3}\right)
\end{align*}
where $\eta=\xi(m'fg_2r_4/m_4'')^2\equiv \overline{\gamma gg_1}m_2(m'fg_2r_4/m_4'')^2 \bmod{r_3}$. Opening the Kloosterman sum we arrive at
\begin{align*}
&\sideset{}{^\star}\sum_{\alpha\bmod{r_3/m_3''}}e\left(\frac{\bar\alpha n}{r_3/m_3''}\right)\\
&\times\sideset{}{^\star}\sum_{\gamma\bmod{r_3}}e\left(\frac{m_3''m'\alpha\gamma gg_1\overline{m_2(m'fg_2r_4/m_4'')^2}}{r_3}-\frac{Mg\gamma d\overline{ph_1r_4}}{r_3}+\frac{\gamma c_2\overline{dmr_4}}{r_3}\right).
\end{align*}
The sum over $\gamma$ is a Ramanujan sum. So we get
\begin{align*}
r_3\sum_{\delta|r_3}\frac{\mu(\delta)}{\delta}
\sideset{}{^\star}\sum_{\substack{\alpha\bmod{r_3/m_3''}\\m_3''m'\alpha gg_1\overline{m_2(m'fg_2r_4/m_4'')^2}\equiv (Mgd^2m- c_2h_1p)\overline{dh_1mpr_4}\bmod{r_3/\delta}}}e\left(\frac{\bar \alpha n}{r_3/m_3''}\right).
\end{align*}
\\

We conclude that the sum in \eqref{41} is dominated by
\begin{align}\label{41+}
\frac{LP}{\tilde{R}^2}\:\frac{dmh(hh_1)^2}{fg_2}&\sum_{r_4|(dmh)^\infty}\:\sum_{m_3''=1}^\infty\:\sum_{m_4''|fg_2r_4m'}\:\mathop{\sum\sum}_{\delta_1,\delta_2=1}^\infty\frac{1}{(\delta_1\delta_2)^2}\sum_{1\leq n<\mathcal{N}}|\lambda(m_3''m_4'',n)|\\
\nonumber \times &\Bigl|\sum_{\substack{P<p<2P\\p\;\text{prime}}}\:\sum_{\substack{|r_3|\asymp \tilde{R}/hh_1\delta_1\delta_2 m_3''r_4\\(gpm_2dmh,r_3)=1}} \:\mathop{\sum}_{\substack{|c_2|\ll \mathcal{C}_2}} \upsilon(\dots)\:\psi_n(\dots)\Bigr|
\end{align}
where the factors $\upsilon(\dots)$ have absolute value smaller than one, and do not depend on $n$, and 
\begin{align*}
\psi_n(\dots)=
\sideset{}{^\star}\sum_{\substack{\alpha\bmod{\delta_2 r_3}\\\delta_1m_3''m'\alpha gg_1\overline{m_2(m'fg_2r_4/m_4'')^2}\equiv (Mgd^2m- c_2h_1p)\overline{dh_1mpr_4}\bmod{r_3m_3''}}}e\left(\frac{\bar \alpha n}{\delta_2 r_3}\right).
\end{align*}
Here $\mathcal{N}$ is as given in \eqref{new-N}, which in terms of the new variables is given by
\begin{align*}
\mathcal{N}=\frac{(fg_2\tilde{R})^3m'M^{36\theta+\varepsilon}}{(hh_1)^3(\delta_1m_3''m_4'')^2L}.
\end{align*}
Note that $\psi_n$ vanishes unless $m_3''|(Mgd^2m- c_2h_1p)$, and in this case we get \begin{align*}
\psi_n(\dots)=
\sideset{}{^\star}\sum_{\substack{\alpha\bmod{\delta_2 r_3}\\\delta_1m' p\eta\alpha \equiv (Mgd^2m- c_2h_1p)/m_3''\bmod{r_3}}}e\left(\frac{\bar \alpha n}{\delta_2 r_3}\right)
\end{align*}
where $\eta=\eta(\dots)=gg_1dh_1mr_4\overline{m_2(m'fg_2r_4/m_4'')^2}$ (which is invertible modulo $\delta_2r_3$).\\

\subsection{Application of Cauchy's inequality and Poisson summation}
Applying the Cauchy inequality and \eqref{ram-on-av}, we see that \eqref{41+} is dominated by
\begin{align}
\label{to-replace}
\sqrt{\frac{L}{\tilde{R}}}M^{48\theta+\varepsilon}P\:\sqrt{hh_1}&\sum_{r_4|(dmh)^\infty}\:\mathop{\sum\sum}_{\delta_1,\delta_2=1}^\infty\;\frac{1}{\delta_1^3\delta_2^2}\:\Bigl\{\sum_{m_3''=1}^{\tilde{R}/\delta_1\delta_2r_4}\:\sum_{m_4''|fg_2r_4m'}\:\Theta_4(\dots)\Bigr\}^{1/2}
\end{align}
where
\begin{align*}
\Theta_4(\dots)=\sum_{1\leq n<\mathcal{N}}\Bigl|\sum_{\substack{P<p<2P\\p\;\text{prime}}}\:\sum_{\substack{|r_3|\asymp \tilde{R}/hh_1\delta_1\delta_2 m_3''r_4\\(gpm_2dmh,r_3)=1}} \:\mathop{\sum}_{\substack{|c_2|\ll \mathcal{C}_2}} \upsilon(\dots)\:\psi_n(\dots)\Bigr|^2.
\end{align*}
Using positivity we now smooth out the $n$-sum and then apply the Poisson summation formula after opening the absolute square. The modulus is $\delta_2r_3r_3'$. We get 
\begin{align*}
\Theta_4(\dots)\leq \mathcal{N}\mathop{\sum\sum}_{\substack{P<p,p'<2P\\p,p'\;\text{prime}}}\:\mathop{\sum\sum}_{\substack{|r_3|, |r_3'|\asymp \tilde{R}/hh_1\delta_1\delta_2 m_3''r_4\\(gpm_2dmh,\: r_3)=1\\(gp'm_2dmh,\: r_3')=1}} \:\mathop{\sum\sum}_{\substack{|c_2|, |c_2'|\ll \mathcal{C}_2}} \upsilon(\dots)\bar{\upsilon}(\dots)\:\sum_{n\in\mathbb{Z}}\mathfrak{C}\mathfrak{I}
\end{align*}
where the character sum (which is again denoted by $\mathfrak{C}$) is given by 
\begin{align}
\label{congruences}
\mathfrak{C}=\mathop{\sideset{}{^\star}\sum_{\alpha\bmod{\delta_2 r_3}}\;\;\sideset{}{^\star}\sum_{\alpha'\bmod{\delta_2 r_3'}}}_{\substack{\delta_1m'p\eta\alpha \equiv (Mgd^2m- c_2h_1p)/m_3''\bmod{r_3}\\\delta_1m'p'\eta\alpha' \equiv (Mgd^2m- c_2'h_1p')/m_3''\bmod{r_3'}\\\bar{\alpha}r_3'- \bar{\alpha}'r_3\equiv n\bmod{\delta_2r_3r_3'}}}1
\end{align}
if $m_3''|(Mgd^2m- c_2h_1p, Mgd^2m- c_2'h_1p')$ (vanishing otherwise) and the integral $\mathfrak{I}$ is the Fourier transform of a smooth bump function. The integral is negligibly small if 
$$
|n|\gg N^\star= \frac{hh_1m_4''^2}{\delta_2(fg_2)^3m'r_4^2}\:\frac{LM^\varepsilon}{\tilde{R}M^{36\theta}}.
$$
We conclude that
\begin{align*}
\Theta_4(\dots)\ll \mathcal{N}\mathfrak{W}+M^{-20130}.
\end{align*}
where
\begin{align}\label{THETA}
\mathfrak{W}=\mathop{\sum\sum}_{\substack{P<p,p'<2P\\p,p'\;\text{prime}}}\:\mathop{\sum\sum}_{\substack{|r_3|, |r_3'|\asymp \tilde{R}/hh_1\delta_1\delta_2 m_3''r_4\\(gpm_2dmh,r_3)=1\\(gp'm_2dmh,r_3')=1}} \:\mathop{\sum\sum}_{\substack{|c_2|, |c_2'|\ll \mathcal{C}_2}}\:\sum_{|n|\ll N^\star} |\mathfrak{C}|.
\end{align}
Recall the definition of $B(\dots)$ from Lemma~\ref{lem-B}. Substituting the above bound into \eqref{to-replace} we conclude the following lemma.\\

\begin{lemma}
\label{lemma2020}
We have
\begin{align*}
B(\dots)\ll \sup\: M^{114\theta+\varepsilon}\tilde{R}P\:&\sum_{\substack{r_4|(dmh)^\infty\\r_4\leq\tilde{R}}}\:\mathop{\sum\sum}_{\delta_1,\delta_2=1}^\infty\;\frac{1}{\delta_1^4\delta_2^2}\:\Bigl\{\sum_{m_3''=1}^{\tilde{R}/\delta_1\delta_2r_4}\:\sum_{m_4''|fg_2r_4m'}\:\frac{\mathfrak{W}}{(m_3''m_4'')^2}\Bigr\}^{1/2}+ M^{-2013}
\end{align*}
where the supremum is taken over all possible vectors $(d,m,h_1,m_2,g,h,g_1,g_2,\beta_1)$ which appear in \eqref{ohhm}, with $m\leq M^{\hat\theta}$, and $dhh_1m_2\leq M^{\theta^\star}$. 
\end{lemma}

\bigskip

\section{A counting problem}

In this section we will estimate the sum that appears in \eqref{THETA}. We are required to count the number of solutions of certain congruence relations. Let us start with the zero frequency $n=0$ contribution which we denote by $\mathfrak{W}_0$. \\

\subsection{The zero frequency}
In this case the last congruence in \eqref{congruences} implies that $r_3=r_3'$ and $\alpha=\alpha'$. The other two congruences now imply that 
$$
p(Mgd^2m- c_2'h_1p')\equiv p'(Mgd^2m- c_2h_1p) \bmod{r_3}.
$$
For any $(p,p',c_2,c_2',r_3)$ satisfying the above congruence we have 
$
\mathfrak{C}\ll \delta_1\delta_2 m', 
$
otherwise $\mathfrak{C}$ vanishes. Now to count the number solutions of the above congruence, we consider two distinct cases. In the first case suppose we have the equality
\begin{align}
\label{equality}
p(Mgd^2m- c_2'h_1p')=p'(Mgd^2m- c_2h_1p).
\end{align}
Then either $p=p'$, in which case we have $c_2=c_2'$, or $p\neq p'$ in which case
$$
Mgd^2m- c_2'h_1p'=p'\ell\;\;\;\text{and}\;\;\;Mgd^2m- c_2h_1p=p\ell
$$
for some $\ell$. The last pair of equalities implies that $pp'|Mgd^2m$ which is ruled out by size considerations and using the fact that $(pp',M)=1$. So the contribution of this case, i.e. when the equality \eqref{equality} holds, to \eqref{THETA} is dominated by 
\begin{align}
\label{zero1}
O\left(\delta_1\delta_2 m'P\mathcal{C}_2\frac{\tilde{R}}{hh_1\delta_1\delta_2m_3''r_4}\right).
\end{align}
Now suppose the equality \eqref{equality} does not hold. Then there are only $O(M^\varepsilon)$ many choices for $r_3$ for any given $(p,p',c_2,c_2')$. So the contribution of this case to \eqref{THETA} is dominated by 
\begin{align}
\label{zero2}
O\left(\delta_1\delta_2m'P^2\mathcal{C}_2^2M^\varepsilon\right).
\end{align}
We will be eventually forced to pick $\theta$ quite small. In particular we will have $\theta<1/450$, so that $P\mathcal{C}_2<\tilde{R}$.\\

\begin{lemma}
\label{lem0}
If $\theta<1/450$ then we have
$$
\mathfrak{W}_0\ll \delta_1\delta_2 \frac{\tilde{R}^2M^{1/2+105\theta+\varepsilon}}{P}.
$$
\end{lemma}

\bigskip

\subsection{Counting the number of solutions of an equation}
Next we consider the contribution of the non-zero frequencies $n\neq 0$, namely
\begin{align*}
\mathfrak{W}-\mathfrak{W}_0=\mathop{\sum\sum}_{\substack{P<p,p'<2P\\p,p'\;\text{prime}}}\:\mathop{\sum\sum}_{\substack{|r_3|, |r_3'|\asymp \tilde{R}/hh_1\delta_1\delta_2 m_3''r_4\\(gpm_2dmh,r_3)=1\\(gp'm_2dmh,r_3')=1}} \:\mathop{\sum\sum}_{\substack{|c_2|, |c_2'|\ll \mathcal{C}_2}}\:\sum_{0<|n|\ll N^\star} |\mathfrak{C}|.
\end{align*} 
In this case we rewrite the congruences in \eqref{congruences} as
\begin{align*}
\delta_1m'p\eta &\equiv \beta\:(Mgd^2m- c_2h_1p)/m_3''\bmod{r_3}\\
\delta_1m'p'\eta &\equiv \beta'\:(Mgd^2m- c_2'h_1p')/m_3''\bmod{r_3'}\\
\beta r_3'- \beta'r_3&\equiv n\bmod{\delta_2r_3r_3'}.
\end{align*} 
Let $u=(r_3,r_3')$ and set $r_3=us$ and $r_3'=us'$ (so $ss'\neq 0$). Then $u|n$ and we write $n=uv$ with $v\neq 0$. It follows that there are at most $(\delta_1\delta_2)^2$ many pairs $(\beta,\beta')$ (so that $\mathfrak{C}\ll (\delta_1\delta_2)^2$) for any given $(c_2,c_2',p,p',s,s',u,v)$ satisfying the congruence conditions
\begin{align}
\label{cong1}\delta_1m'm_3''p\eta s'&\equiv v\:(Mgd^2m- c_2h_1p)\bmod{s}\\
\label{cong2}\delta_1m'm_3''p'\eta s &\equiv -v\:(Mgd^2m- c_2'h_1p')\bmod{s'}
\end{align}
and
\begin{align*}
&\delta_1m'm_3''p\eta(Mgd^2m- c_2'h_1p') s'- \delta_1m'm_3''p'\eta(Mgd^2m- c_2h_1p) s\\
&\equiv v(Mgd^2m- c_2h_1p)(Mgd^2m- c_2'h_1p')\bmod{u}.
\end{align*} 
Recall that $\eta=\eta_1\bar{\eta}_2$ where $\eta_1=dgg_1h_1mm_4''^2$ and $\eta_2=(fg_2m')^2m_2r_4$. Let $\mathfrak{V}$ be the number of solutions $(c_2,c_2',p,p',s,s',u,v)$ of the above set of three congruences within the desired range. Then
$$
\mathfrak{W}-\mathfrak{W}_0\ll (\delta_1\delta_2)^2\:\mathfrak{V}.
$$\\

We will first tackle a degenerate case. We will denote the contribution of this part to $\mathfrak{V}$ by $\mathfrak{V}_{e}$. Suppose we have the equality
\begin{align}
\label{equation-of-variety}
&\delta_1m'm_3''p\eta_1(Mgd^2m- c_2'h_1p') s'- \delta_1m'm_3''p'\eta_1(Mgd^2m- c_2h_1p) s\\
\nonumber &= v\eta_2(Mgd^2m- c_2h_1p)(Mgd^2m- c_2'h_1p').
\end{align}
Then all the three congruences above hold for any $u$. To count the number of solutions of this equation, suppose we are given $(p,c_2,s)$. Then one may determine $(p',c_2')$ from the congruence
\begin{align*}
\delta_1m'm_3''p'\eta_1(Mgd^2m- c_2h_1p) s\equiv 0\bmod{(Mgd^2m- c_2'h_1p')}.
\end{align*}
This boils down to the divisibility condition 
\begin{align*}
(Mgd^2m- c_2'h_1p')\:|\:\delta_1m'm_3''p'\eta_1(Mgd^2m- c_2h_1p) s.
\end{align*}
Note that the right hand side is non-zero as $p\nmid Mgd^2m$ ($p\nmid M$ by our choice of $p$, and $p\nmid gdm$ due to size restrictions). Hence there are $O(M^\varepsilon)$ many choices for the pairs $(c_2',p')$. It remains to estimate the number of $s'$ and $v$. Suppose $(s',v)$ and $(\bar{s}',\bar{v})$ are two such pairs with $v\neq \bar{v}$. Then we have
\begin{align*}
&\delta_1m'm_3''p\eta_1(Mgd^2m- c_2'h_1p') (s'-\bar{s}')\\
&= (v-\bar{v})\eta_2(Mgd^2m- c_2h_1p)(Mgd^2m- c_2'h_1p'),
\end{align*}
which reduces to  
\begin{align*}
\delta_1m'm_3''p\eta_1 (s'-\bar{s}')= (v-\bar{v})\eta_2(Mgd^2m- c_2h_1p),
\end{align*}
and consequently 
\begin{align*}
(Mgd^2m- c_2h_1p)\:|\:\delta_1m'm_3''\eta_1 (s'-\bar{s}').
\end{align*}
Given $s'$, there are 
$$
O\left(1+\frac{\tilde{R}(Mgd^2m- c_2h_1p,\delta_1m'm_3''\eta_1)}{uhh_1\delta_1\delta_2m_3''r_4(Mgd^2m- c_2h_1p)}\right)
$$
many choices for $\bar{s}'$. It follows that the number of solutions of the equation \eqref{equation-of-variety} is bounded by 
\begin{align}
\label{first-second-term}
M^\varepsilon\sum_{1\leq u \ll N^\star}&\mathop{\sum}_{\substack{P<p<2P\\p\;\text{prime}}}\:\mathop{\sum}_{\substack{|s|\asymp \tilde{R}/uhh_1\delta_1\delta_2 m_3''r_4}}\\
\nonumber & \:\mathop{\sum}_{\substack{|c_2|\ll \mathcal{C}_2}}\left(1+\frac{\tilde{R}(Mgd^2m- c_2h_1p,\delta_1m'm_3''\eta_1)}{uhh_1\delta_1\delta_2m_3''r_4(Mgd^2m- c_2h_1p)}\right).
\end{align}
The contribution of the first term is given by
\begin{align*}
M^\varepsilon\sum_{1\leq u \ll N^\star}\mathop{\sum}_{\substack{P<p<2P\\p\;\text{prime}}}\:\mathop{\sum}_{\substack{|s|\asymp \tilde{R}/uhh_1\delta_1\delta_2 m_3''r_4}} \:\mathop{\sum}_{\substack{|c_2|\ll \mathcal{C}_2}}1\ll M^\varepsilon\frac{P\tilde{R}\mathcal{C}_2}{hh_1\delta_1\delta_2m_3''r_4}.
\end{align*}
On the other hand the second term gives
\begin{align*}
M^\varepsilon\sum_{1\leq u \ll N^\star}\left(\frac{\tilde{R}}{uhh_1\delta_1\delta_2m_3''r_4}\right)^2\mathop{\sum}_{\substack{P<p<2P\\p\;\text{prime}}}\:\mathop{\sum}_{\substack{|c_2|\ll \mathcal{C}_2}}\frac{(Mgd^2m- c_2h_1p,\delta_1m'm_3''\eta_1)}{(Mgd^2m- c_2h_1p)}.
\end{align*}
Consider the inner sums over $c_2$ and $p$. The part of the sum with $c_2=0$ boils down to
\begin{align*}
\mathop{\sum}_{\substack{P<p<2P\\p\;\text{prime}}}\:\frac{(Mgd^2m,\delta_1m'm_3''\eta_1)}{Mgd^2m}\ll\frac{P}{M}.
\end{align*}
The last inequality follows as $(M,\delta_1m'm_3''\eta_1)=1$ (recall that we are assuming that $M$ is a prime number). For $c_2\neq 0$, we set $j=Mgd^2m-c_2h_1p$, which ranges over a set of non-zero integers. For any given $j$ there are at most $O(M^\varepsilon)$ many pairs $(p,c_2)$ (as $c_2\neq 0$) such that $j=Mgd^2m-c_2h_1p$. So the sum is dominated by
$$
M^\varepsilon\sum_{0<|j|<J}\frac{(j,\delta_1m'm_3''\eta_1)}{|j|}\ll M^\varepsilon\sum_{\substack{\zeta|\delta_1m'm_3''\eta_1\\\zeta<J}}\;\sum_{0<|j|<J/\zeta}\frac{1}{|j|}\ll M^\varepsilon
$$ 
(where $J=M^{2013}$ say). Since $\sqrt{M}<P<M$ we see that the expression in \eqref{first-second-term} is dominated by
\begin{align}
\label{to-compare}
O\left(\frac{\tilde{R}^2M^\varepsilon}{(hh_1\delta_1\delta_2m_3''r_4)^2}\left(1+\frac{hh_1\delta_1\delta_2m_3''r_4 P\mathcal{C}_2}{\tilde{R}}\right)\right)=O\left(\frac{\tilde{R}^2}{\delta_1\delta_2}M^\varepsilon\right).
\end{align}
In the last equality we again assumed that $\theta<1/450$. This is the bound for the number of solutions $\mathfrak{V}_e$ of the equation \eqref{equation-of-variety}. \\

\begin{lemma}
\label{lem-eqn}
For $\theta<1/450$ we have
$$
(\delta_1\delta_2)^2\:\mathfrak{V}_e \ll \delta_1\delta_2\:\tilde{R}^2\:M^\varepsilon.
$$
\end{lemma}

\bigskip

\subsection{The generic count}
Now we will count the number of solutions of the pair of congruences \eqref{cong1} and \eqref{cong2}, which are not coming from the equation \eqref{equation-of-variety}. We write the first congruence \eqref{cong1} as an equality
\begin{align}
\label{the-eqn}
\delta_1m'm_3''\eta_1\: ps'= v\:(Mgd^2m- c_2h_1p)\eta_2 + es.
\end{align} 
From the second congruence  \eqref{cong2} now it follows that
\begin{align*}
&(Mgd^2m- c_2h_1p)\delta_1m'm_3''p'\eta_1 s\\
 &\equiv -(Mgd^2m- c_2'h_1p')(\delta_1m'm_3''\eta_1\: ps'-es)\bmod{s'},
\end{align*} 
and hence
\begin{align*}
s'|(Mgd^2m- c_2h_1p)\delta_1m'm_3''p'\eta_1 - (Mgd^2m- c_2'h_1p')e.
\end{align*} 
So once $p$, $p'$, $c_2$, $c_2'$ and $e$ are given there are $O(M^\varepsilon)$ many choices for $s'$, unless 
\begin{align*}
(Mgd^2m- c_2h_1p)\delta_1m'm_3''p'\eta_1= (Mgd^2m- c_2'h_1p')e.
\end{align*} 
Now suppose we have obtained $p$, $p'$, $c_2$, $c_2'$, $e$ and $s'$. Next we count the number of $v$ satisfying the congruence
\begin{align*}
\delta_1m'm_3''\eta_1\: ps'\equiv v\:(Mgd^2m- c_2h_1p)\eta_2 \bmod{e}.
\end{align*} 
This is bounded by 
\begin{align}
\label{count-for-v}
O\left(1+\frac{N^\star(Mgd^2m- c_2h_1p,e)}{|e|}\right).
\end{align}
Finally $s$ is uniquely determined by the equation \eqref{the-eqn} if $e\neq 0$, and there are $O(M^\varepsilon)$ many choices for $u$ as we are only interested in solutions not satisfying the equation \eqref{equation-of-variety}. Observe that
\begin{align}
\label{last-to-name}
\mathop{\mathop{\sum\sum}_{p,p'\sim P}\mathop{\sum\sum}_{|c_2|,|c_2'|\ll \mathcal{C}_2}\sum_{0<|e|\ll E}}_{(Mgd^2m- c_2h_1p)\delta_1m'm_3''p'\eta_1\neq (Mgd^2m- c_2'h_1p')e}&\left(1+\frac{N^\star(Mgd^2m- c_2h_1p,e)}{|e|}\right)\\
\nonumber &\ll M^\varepsilon P^2\mathcal{C}_2^2(E+N^\star).
\end{align}
Considering the sizes of the variables in the equation \eqref{the-eqn} it follows that
$$
E\ll \delta_1dgg_1hh_1mm'm_3''m_4''^2\:PM^{\varepsilon}.
$$
Also we have
$$
N^\star\ll \frac{hh_1m_4''^2}{\delta_2(fg_2)^3m'r_4^2}\:\frac{P^2M^\varepsilon}{M^{1/2+19\theta}},
$$
and consequently (as $P>\sqrt{N}$ and $\theta<1/450$)
$$
(E+N^\star)\ll \delta_1dgg_1hh_1mm'm_3''m_4''^2\:P^2M^{-1/2-19\theta+\varepsilon}.
$$
Hence the contribution of these non-zero frequencies to $\mathfrak{V}$ is given by
$$
\mathfrak{V}_g\ll \delta_1dgg_1hh_1mm'm_3''m_4''^2\:P^2\mathcal{C}_2^2 P^2M^{-1/2-19\theta+\varepsilon}
$$\\

\begin{lemma}
\label{lem-gen}
For $\theta<1/450$ we have
$$
(\delta_1\delta_2)^2\mathfrak{V}_g \ll \delta_1^3\delta_2^2m_3''m_4''^2\:\tilde{R}^2M^{1/2+103\theta+\varepsilon}.
$$
\end{lemma}

\bigskip
\subsection{Two degenerate cases}
Now we are left with two degenerate cases whose contributions are not included in $\mathfrak{V}_e+\mathfrak{V}_g$. First suppose we have $p$, $p'$, $c_2$, $c_2'$ and $e$ such that
$$
(Mgd^2m- c_2h_1p)\delta_1m'm_3''p'\eta_1= (Mgd^2m- c_2'h_1p')e.
$$
Suppose we are given the pair $(p,c_2)$. Then we can get at most $O(M^\varepsilon)$ many pairs $(p',c_2')$ satisfying the divisibility condition
$$
(Mgd^2m- c_2'h_1p')\:|\:(Mgd^2m- c_2h_1p)\delta_1m'm_3''\eta_1.
$$
Then once we have obtained $(p,p',c_2,c_2')$ there are $O(M^\varepsilon)$ many choices for $e$. Next we count the number of solutions $(s,s',v)$ of the linear equation \eqref{the-eqn}. Viewing it as a congruence modulo $p$ we see that the number of solutions is 
$$
O\left(M^\varepsilon\frac{N^\star}{u}\left(1+\frac{\tilde{R}}{Puhh_1\delta_1\delta_2 m_3''r_4}\right)\right).
$$
As above we now have \eqref{count-for-v} many choices for $v$ and $O(M^\varepsilon)$ many choices for $s$. So the contribution of this case to $\mathfrak{V}$ is bounded by
\begin{align*}
&M^\varepsilon\sum_{1\leq u\ll N^\star}\mathop{\sum}_{p\sim P}\mathop{\sum}_{|c_2|\ll \mathcal{C}_2}\:\frac{N^\star}{u}\left(1+\frac{\tilde{R}}{Puhh_1\delta_1\delta_2 m_3''r_4}\right)\\
&\ll M^\varepsilon \:\left(P\mathcal{C}_2N^\star+\frac{\mathcal{C}_2N^\star\tilde{R}}{hh_1\delta_1\delta_2m_3''r_4}\right).
\end{align*}
Now
$$
\frac{\tilde{R}}{P^2\mathcal{C}_2hh_1\delta_1\delta_2m_3''r_4}\ll \frac{1}{\sqrt{M}}.
$$
We see that the contribution of this degenerate case is dominated by \eqref{last-to-name}.\\

Finally we consider the case $e=0$. So that the congruence \eqref{cong1} reduces to an equation. Using symmetry we reduce the problem to counting the number of solutions of the pair of equations
\begin{align*}
\delta_1m'm_3''p\eta s'&=v\:(Mgd^2m- c_2h_1p)\\
\delta_1m'm_3''p'\eta s &= -v\:(Mgd^2m- c_2'h_1p').
\end{align*}
Once a (non-zero) value for $s'$ is given, we have $O(M^\varepsilon)$ many choices for $(p,c_2)$ because of the divisibility condition 
\begin{align*}
(Mgd^2m- c_2h_1p)\:|\:\delta_1m'm_3''\eta s'.
\end{align*}
Then $v$ is determined from the first equation. Similarly $s$ determines the other variables. So the contribution of this part is dominated by
$$ 
M^\varepsilon \:\frac{\tilde{R}^2}{(hh_1\delta_1\delta_2m_3''r_4)^2}.
$$
Compare with \eqref{to-compare}. So this is also satisfactory.\\

Finally we observe that the bounds from Lemma~\ref{lem0} and Lemma~\ref{lem-eqn} can be absorbed in the bound given in Lemma~\ref{lem-gen}. So we conclude the following lemma.\\

\begin{lemma}
\label{lemfinalfinal}
For $\theta<1/450$ we have
$$
\mathfrak{W} \ll \delta_1^3\delta_2^2m_3''m_4''^2\:\tilde{R}^2M^{1/2+103\theta+\varepsilon}.
$$
\end{lemma}

\bigskip 

\section{Conclusion}

Substituting the bounds from Lemma~\ref{lemfinalfinal} into Lemma~\ref{lemma2020} we get 
$$
B(\dots)\ll \tilde{R}^2PM^{1/4+166\theta+\varepsilon}.
$$
Then from Lemma~\ref{lem-B} we conclude that
for $\sqrt{N}<P<M$ and $\theta<1/450$ we have
\begin{align*}
\sup_{\mathfrak{m}'\mathfrak{d}\leq M^{\theta^\star}}\:|\mathcal{O}(C,m;L,\mathfrak{m}',\mathfrak{d})| \ll \sqrt{N}\:\frac{M^{3/2+201\theta+\varepsilon}}{m^2P}.
\end{align*} 
From Corollary~\ref{cor-lemma15} it now follows that
\begin{align*}
\mathcal{O}^\star(C,m)\ll M^\varepsilon \:\frac{\sqrt{N}}{m^{3/2}}\left(\frac{M^{3/2+201\theta}}{P}+ M^{3/4-\theta/2}\right)
\end{align*}
So the optimal choice for $P$ is given by
$$
P=M^{3/4+403\theta/2}.
$$
This is allowed if $M^{3/4+403\theta/2}<M$ or $\theta<1/806$. Consequently we get 
\begin{align*}
\mathcal{O}^\star(C,m)\ll \frac{\sqrt{N}}{m^{3/2}}M^{3/4-1/1612+\varepsilon}.
\end{align*}
Substituting into Corollary~\ref{cor-previous-1115} and summing over $m$, we get
$$
\mathcal{F}\ll \sqrt{N}M^{3/4-1/1612+\varepsilon}.
$$
From Lemma~\ref{lem-31}, we see that the same bound holds for $\mathcal{O}$ as well. Consequently
$$
S^\star(N)\ll \sqrt{N}M^{3/4-1/1612+\varepsilon}.
$$
Substituting this into Corollary~\ref{first-cor}, we obtain the theorem.

\end{document}